%% file: ADMM_Tan_HAL3.tex
\documentclass[a4paper,10pt]{article}

\input{preambuleGeneral}
\newcommand{\titre}{Linear Convergence Rates for Variants of the ADMM in Smooth Cases}

\definecolor{couleur}{RGB}{109,76,90} 
\input{preambuleMath}

\input{preambuleMiseenpage}
\begin{document}

\begin{montitre}
Acceleration of saddle-point methods

in smooth cases
\end{montitre}
\vspace{-0.75cm}

\begin{center}
Pauline \textsc{Tan}

CMAP, \'Ecole polytechnique, CNRS,

Université Paris-Saclay, 91128, Palaiseau, France
\end{center}
\bigskip
\renewcommand{\abstractname}{\textcolor{couleur}{Abstract}}

\begin{abstract}
  In the present paper we propose a novel convergence analysis of the Alternating Direction Methods of Multipliers (ADMM), based on its equivalence with the overrelaxed Primal-Dual Hybrid Gradient (oPDHG) algorithm. We consider the smooth case, which correspond to the cas where the objective function can be decomposed into one differentiable with Lipschitz continuous gradient part and one strongly convex part. An accelerated variant of the ADMM is also proposed, which is shown to converge linearly with same rate as the oPDHG.
\end{abstract}

\section{Introduction}
\subsection{Context}
The Alternating Direction Methods of Multipliers (ADMM) is a widely-used method aimed at minimizing constrained problems of form
\begin{equation}
\min_{\substack{(x,z)\in X\times Z\\Ax + Bz = c}} g(x) + h(z).
\label{eq:Pprimalcontraint}
\end{equation}
The objective function is separable in~$(x,z)$ with~$g:X\to\mathbb{R}\cup\{+\infty\}$ and~$h:Z\to\mathbb{R}\cup\{+\infty\}$ two closed convex functions.
The constraint involves two linear operators~$A:X\to Y$ and~$B:Z\to Y$ and a constant~$c\in Y$. In this work, $X$, $Z$, and~$Y$ are finite-dimensional real Hilbert spaces.
The ADMM was initially introduced in the mid-70's by \textsc{Gabay-Mercier}~\cite{gabay1976dual} and by \textsc{Glowinski-Marrocco}~\cite{glowinski1975approximation}. It considers the augmented Lagrangian associated to problem~\eqref{eq:Pprimalcontraint}
\begin{equation}
L_{\tau}(x,z;y) := g(x) + h(z) + \langle Ax+Bz-c,y \rangle + \frac{1}{2\,\tau}\,\lVert Ax+Bz-c\rVert^2
\label{eq:lagrangien}
\end{equation}
for~$\tau>0$ which leads to solve the saddle-point problem
\begin{equation}
\min_{(x,z)\in X\times Z}\sup_{y\in Y}\quad L_{\tau}(x,z;y)
\label{eq:ADMM}
\end{equation}
instead of the initial problem.
One particular instance of these so-called augmented Lagrangian methods uses \textsc{Uzawa}'s method to solve~\eqref{eq:ADMM}. Namely, the \emph{method of multipliers} tackles this problem by alternating an exact minimization on the primal variable~$(x,z)$ and a gradient ascent step on the dual variable~$y$. In such a method, the minimization step couples the primal variables. 
To decouple them, one may consider splitting this step into two partial minimizations, one over~$x$ and another over~$z$. These two minimization can be done simultaneously, from the same initial points, or, in the case of the ADMM, one after the other, with an update in between. This leads to the following algorithm
\begin{equation}
\begin{cases}
x_{n+1} = \displaystyle\arg\min_{x\in X} L_{\tau}(x,z_n;y_n)\\[3mm]
z_{n+1} = \displaystyle\arg\min_{z\in Z} L_{\tau}(x_{n+1},z;y_n)\\[3mm]
y_{n+1} = \displaystyle y_n + \frac{1}{\tau}\,(Ax_{n+1} - z_{n+1}).
\end{cases}
\label{eq:algoADMM0}
\end{equation}

This method can be proved to be linked to another famous algorithm, which is known as the Primal-Dual Hybrid Gradient (PDHG) method~\cite{zhu2008efficient}. 
The PDHG method tackles saddle-point problems by alternating gradient descent steps and gradient ascent steps. Such problems arise while considering a primal-dual formulation of a convex minimization problem, in a splitting strategy for instance. 
A noteworthy feature of the PDHG method is that it can be accelerated thanks to an overrelaxation step \emph{\`a la} \textsc{Nesterov}~\cite{nesterov1983method} on one of the variables~\cite{pock2009algorithm, chambolle2011first,esser2010general,chambolle2015ergodic}, which leads to the overrelaxed~PDHG (oPDHG).

The ADMM has been intensively studied in the past years.
One may see for instance a comprehensive review in~\cite{boyd2011distributed}. 
The key point is the convergence of the algorithm and its convergence rate. 
Under assumptions on the matrix ranks and / or the regularity of the objective functions~$g$ and~$h$, linear rates can be achieved~\cite{hong2012linear}.
Eventually, some accelerated variants of the ADMM have been proposed~\cite{deng2016global,davis2014faster}.

As a recent developpement, we should mention~\cite{he2014convergence}, which also studied the convergence of the PDHG method and derived optimal step size choice, when only one function assumed to be strongly convex.

\subsection{Contribution of this paper}


In this paper, we provide a new analysis of the ADMM based on the equivalence between the ADMM and the oPDHG method. More specifically, we use the analysis to derive convergence rate for the~ADMM in a case we refer to be \emph{smooth}. We indeed made restrictive assumptions on the initial problem~\eqref{eq:Pprimalcontraint}, which implies that we consider the following particular instance of~\eqref{eq:Pprimalcontraint}:
\begin{equation}
\min_{\substack{(x,z)\in X\times Y\\Ax = z}}  g(x)+h(z).
\label{eq:Pprimalcontraint2}
\end{equation}
which may be rewritten as the unconstrained composite problem
\begin{equation}
\min_{x\in X}  g(x)+h(Ax)
\label{eq:Pprimal0}
\end{equation}
with regularity assumptions on $g$, which is supposed to be strongly convex, and $h$, which has a {Lipschitz} gradient. We first establish new linear ergodic convergence rates of the oPDHG by generalizing the proofs of \cite{chambolle2011first,chambolle2015ergodic}. This leads to a linear rate for the ADMM under these assumptions.
Then, we introduce a slight variant of the~ADMM which leads to a better rate, by relaxing the choice of the parameters in the convergence proof of the oPDHG method.

The reason why we only consider the case~$B=-\text{Id}$ and~$c=0$ is that, otherwise, as the map~$y\mapsto h^*(B^*y)$ will be supposed to be strongly convex, this implies that~$\nabla f$ is Lipschitz continuous and that~$B$ is invertible. Such conditions are artificial when~$B$ is not~$-\text{Id}$. However, the interested reader will easily extend our result to this case. 
Moreover, problems of standard form~\eqref{eq:Pprimal0} often arise in many contexts, and thus can justify a special study by themselves.

\subsection{Structure of the paper}
This paper is organized as follows. In Section~\ref{sec:eq}, we recall the equivalence between the ADMM and the oPDHG method. We also define what we call the \emph{smooth case}, which is the case we will consider throughout this paper. In Section~\ref{sec:cvPDHG}, we establish two linear convergence results for the oPDHG, and we provide the best parameter choice in the case where the overrelaxation parameter is fixed to be $1$ or left unconstrained. In Section~\ref{sec:cvADMM}, we exploit the equivalence between the ADMM and the oPDHG to derive from the results of Section~\ref{sec:cvPDHG} new linear convergence rate for the ADMM. We also propose a slight variant of the ADMM, which leads in the best case to the same convergence rate as the oPDHG method. In Section~\ref{sec:relations}, we compare our results with some found in the literature for the classical ADMM or variants, in the case where the assumptions made on the problem yield a linear convergence rate. Those assumptions do not necessary include the smooth case studied here. Eventually, in Section~\ref{sec:app}, we applied our accelerated ADMM on two problems, and compared its convergent with the unaccelerated ADMM, the oPDHG and an adaptation of \textsc{Beck} and \textsc{Teboulle}'s FISTA~\cite{beck2009fast} for the strongly convex case~\cite{nesterov2004introductory,chambolle2016introduction}.

\section{Equivalence between the ADMM and the oPDHG}\label{sec:eq}
\subsection{Initial primal problem}

Let~$X$ and~$Y$ be two finite-dimensional real Hilbert spaces. The inner product is denoted by~$\langle\cdot,\cdot\rangle$ and~$\lVert\cdot\rVert$ stands for the induced norm.
We recall that we consider the minimization problem
\begin{equation}
\min_{x\in X}  \Big\{f(x) := g(x)+h(Ax)\Big\}
\label{eq:Pprimal1}
\end{equation}
where~$g:X\to\mathbb{R}\cup\{+\infty\}$ and~$h:Y\to\mathbb{R}\cup\{+\infty\}$ are proper, convex, and lower semi-continuous (l.s.c.) functions. 
The map~$A:X\to Y$ is a continuous linear operator. Its adjoint is denoted by~$A^*$ and it is supposed to be bounded, of norm~$L_A$
\begin{equation}
L_A := \lVert A\rVert = \sup_{x\in X, \lVert x\rVert \leq 1} \lVert Ax\rVert.
\label{eq:normeA}
\end{equation}


\subsection{Equivalence with oPDHG}\label{sub:equivalence}
Let us briefly recall how the ADMM is connected to the oPDHG algorithm, by rewritting the ADMM iterations~\eqref{eq:algoADMM0} applied on Problem \eqref{eq:Pprimal1}. Ignoring the constant terms in the minimization steps, we obtain
\begin{equation}
\begin{cases}
x_{n+1} = \displaystyle\arg\min_{x\in X}\left\{g(x) + \langle Ax,y_n\rangle + \frac{1}{2\,\tau}\,\lVert Ax-z_n\rVert^2 \right\}\\[4mm]
z_{n+1} = \displaystyle\arg\min_{z\in Y}\left\{h(z) - \langle z,y_n\rangle + \frac{1}{2\,\tau}\,\lVert Ax_{n+1}-z\rVert^2 \right\}\\[5mm]
y_{n+1} = \displaystyle y_n + \frac{1}{\tau}\,(Ax_{n+1} - z_{n+1}).
\end{cases}
\label{eq:algoADMM}
\end{equation}
Defining~$\xi_{n+1}:=Ax_{n+1}$ and introducing the map
\begin{equation}
g_A(\xi):=\inf_{x\in X,Ax=\xi} g(x)
\end{equation}
we can make a change of variable in the~$x$-update and rewrite the updates of~$x_{n+1}$ and~$y_{n+1}$ thanks to proximity operators. This yields
\begin{equation}
\begin{cases}
\xi_{n+1} = \text{prox}_{\tau g_A}\big(\xi_n - \tau\,\bar y_n\big)\\[2mm]
y_{n+1} = \text{prox}_{h^*/\tau}(y_n+\xi_{n+1}/\tau)\\[2mm]
\bar y_{n+1} = y_{n+1}+(y_{n+1}-y_n)
\end{cases}
\label{eq:algoADMM-PD}
\end{equation}
and the $z$-update is given by $z_{n+1} = \xi_{n+1} - \tau\,(y_{n+1}-y_n)$.
This primal-dual algorithm has been studied in~\cite{chambolle2011first}. It can be interpreted as an PDHG algorithm with an additional overrelaxation step (of parameter~$1$) on the dual variable. It solves the saddle-point problem
\begin{equation}
\min_{\xi\in Y}\sup_{y\in Y} \Big\{ g_A(\xi) + \langle \xi,y \rangle - h^*(y) \Big\}
\label{eq:Pprimaldual0}
\end{equation}
which is of general form 
\begin{equation}
\min_{\xi\in Z}\sup_{y\in Y} \Big\{ \mathcal{L}(\xi;y):= G(\xi) + \langle K\xi,y \rangle - H^*(y) \Big\}
\label{eq:Pprimaldual}
\end{equation}
with~$Z=Y$, $K=\text{Id}$,~$G=g_A$ and~$H=h$. Note that~\eqref{eq:Pprimaldual} is the primal-dual formulation of the minimization problem
\begin{equation}
\min_{\xi\in Z} \Big\{ G(\xi) + H(K\xi) \Big\}.
\label{eq:Pgeneral}
\end{equation}

\subsection{Smooth case}\label{sub:assumptions}
From now on, we consider the \emph{smooth case}. 
In the initial primal problem~\eqref{eq:Pprimal1}, the functions~$g$ and~$h^*$ are both supposed to be strongly convex, with respective parameter~$\gamma>0$ and~$\delta>0$. 
We recall that a function~$f:X\to\mathbb{R}\cup\{+\infty\}$ is strongly convex of parameter~$\alpha>0$ ($f$ is also said to be~$\alpha$-convex) if for any~$x_1,x_2\in X$ and~$p\in\partial f(x_1)$
\begin{equation}
f(x_2)\geq f(x_1) + \langle p,x_2-x_1\rangle + \frac{\alpha}{2}\,\lVert x_2-x_1\rVert^2
\label{eq:fortconvex}
\end{equation}
where~$\partial f(x_1)$ denotes the subdifferential of~$f$ at point~$x_1$. 
One can easily check that if~$f$ is $\alpha$-convex, then its convex conjugate $f^*$ is differentiable, with a Lipschitz continuous gradient, of constant $1/\alpha$. 

Let us study the regularity of Problem \eqref{eq:Pprimaldual}. The assumptions made above imply obviously that~$H^*$ is $\delta$-convex.
Moreover, it is easy to show that~$G^*$ is differentiable and that~$\nabla G^*$ is Lipschitz continuous with constant~$L_A^2/\gamma$, which follows from
\begin{equation}
g_A^*(y+t)
= g^*\big(A^*(y+t)\big)
= g^*(A^*y) + \langle \nabla g^*(A^*y), A^*t\rangle + o(\lVert A^*t\rVert)
\end{equation}
since~$g$ is~$\gamma$-convex. Hence,~$G$ is~$\gamma/L_A^2$-convex. Let~$\tilde\gamma = \gamma/L_A^2$ and~$\tilde\delta =\delta$. 

We define $\kappa_f:=L_A^2/(\delta\gamma)$ the \emph{condition number} of $f$ as the ratio between $L_A^2/\delta$ the Lipschitz constant of the smooth part $h(K\cdot)$ and $\gamma$ the strong convexity parameter of the non-smooth part $g$. In the case where $f$ is both smooth with $\nabla f$ lipschitz continuous and strongly convex, this definition recovers the one usually used in such cases and the condition number is always larger than $1$. In the general case, it can be less than $1$. When $\kappa_f$ is large, the function is said ill-conditioned.

\subsection{Forward-backward splitting}
If $h$ is differentiable, it is possible to consider a forward-backward splitting (FBS) strategy to solve problem~\eqref{eq:Pprimal1}. The FBS applied on the sum $f=g+h(A\cdot)$ gives updates of form
\begin{equation}
x_{n+1} = \text{prox}_{\tau g}\big(x_n - \tau\,A^*\,\nabla h(Ax_n)\big).
\label{eq:updateFBS}
\end{equation}
Hence, choosing to use the FBS instead of the ADMM or the oPDHG method suggests that $\nabla h$ is supposed to be easier to compute than $\text{prox}_h$.

A variant of the FBS is FISTA~\cite{beck2009fast}, which adds an extra overrelaxation step. It can be adapted to solve for strongly convex problems following~\cite{nesterov2004introductory}, see~\cite[Appendix B]{chambolle2016introduction} for details. In other terms, the updates~\eqref{eq:updateFBS} are replaced by
\begin{equation}
\begin{cases}
x_{n+1} = \text{prox}_{\tau g}\big(\bar x_n - \tau\,A^*\,\nabla h(A\bar x_n)\big)\\
\bar x_{n+1} = x_{n+1} + \theta_{n+1}\,(x_{n+1}-x_n)
\end{cases}
\label{eq:updateFBS2}
\end{equation}
where the variable overrelaxation parameter $\theta_n$ is chosen in the strongly convex case by letting
\begin{equation}
t_{n+1} = \frac{1-q\,t_n^2+\sqrt{(1-q\,{t_n}^2)^2+4\,{t_n}^2}}{2}
\end{equation}
for $q=\tau\gamma/(1+\tau\gamma)$ for $\tau\in (0,\delta/L_A^2]$.  We will refer to this algorithm as `strongly convex FISTA' in this paper. Then,
\begin{equation}
\theta_n = \big(1+\tau\gamma(1-t_{n+1})\big)\,\frac{t_n-1}{t_{n+1}}.
\label{eq:thetaFISTA}
\end{equation}
In the non-strongly convex case ($\gamma=\delta = 0$), the quantity $q$ is null, and the resulting updates of $t_n$ and $\theta_n$ are those in the original paper of \textsc{Beck} and \textsc{Teboulle}. When~$g$ is assumed to be strongly convex and $h(A\cdotp)$ has a $L^2_A/\delta$-\textsc{Lipschitz} gradient, the convergence rate for the objective error of this algorithm has been proved to be linear. In the case where $t$ (and thus, $\theta$) is chosen to be constant
\begin{equation}
t_n = t = \frac{1}{\sqrt{q}}\qquad\text{and}\qquad
\theta_n = \theta = (1-\sqrt{q})^2 \,\frac{1+\tau\,\gamma}{1-\tau\,\gamma}
\label{eq:FISTApasconstant}
\end{equation}
then the linear rate is of parameter $1-\sqrt{q}$~\cite[Remark B.2]{chambolle2016introduction}. This rate is minimal when~$\tau$ is maximal and equals
\begin{equation}
\omega
 = 1 - \sqrt{\frac{\delta\gamma/L_A^2}{1+\delta\gamma/L_A^2}}
 = 1 - \sqrt{\frac{1}{\kappa_f+1}}.
\end{equation}
Note that if $g$ is $\gamma$-convex then this so is $f$. Hence, the optimality condition on $x^*$ coupled with the strong convexity inequality recalled in~\eqref{eq:fortconvex} yields
\begin{equation}
f(x_n)-f(x^*) \geq \frac{\gamma}{2}\,\lVert x_n - x^*\rVert^2
\end{equation}
that is, a linear convergence for the objective error implies a linear convergence of at least same rate for the convergence of the primal iterate $x_n$.

\section{Convergence of oPDHG in the smooth case}\label{sec:cvPDHG}
In this section, we establish the general convergence proof of the following algorithm
\begin{equation}
\begin{cases}
y_{n+1} = \text{prox}_{\sigma H^*}(y_n+\sigma\,K\bar\xi_n)\\[2mm]
\xi_{n+1} = \text{prox}_{\tau G}\big(\xi_n - \tau\,K^*y_{n+1}\big)\\[2mm]
\bar \xi_{n+1} = \xi_{n+1}+\theta\,(\xi_{n+1}-\xi_n).
\end{cases}
\label{eq:algoADMM-PD22}
\end{equation}
which aims at solving problem~\eqref{eq:Pprimaldual}, in the general case where~$K:Z\to Y$ is bounded of norm~$L_K$,~$G:Z\to\mathbb{R}\cup\{+\infty\}$ is~$\tilde\gamma$-convex and~$H^*:Y\to\mathbb{R}\cup\{+\infty\}$ is~$\tilde\delta$-convex. The step sizes~$\tau,\sigma>0$ and the relaxation parameter~$0<\theta\leq 1$ are to be specified.

When~$\theta=0$, this algorithm is known as the PDHG method~\cite{zhu2008efficient}. It consists in a proximal gradient ascent step for the dual variable, followed by a proximal gradient descent step for the primal variable. 
The overrelaxation step has been added in~\cite{pock2009algorithm} for minimizing the \textsc{Mumford-Shah} functional, and studied in a wider framework in~\cite{esser2010general} and more recently in~\cite{chambolle2015ergodic}.
The case~$\theta=1$ and~$\tau\sigma=1$ corresponds to the equivalence with the ADMM, as recalled in the previous section.
When~$\theta=1$ and~$\tau\sigma\neq1$, the iterations are equivalent to the ADMM with an additional proximal term~\cite{chambolle2011first}, which leads to a preconditioned version of the ADMM~\cite{esser2010general}.
 
Now we can formulate our main result.

\begin{theorem} \label{thm1}
Assume problem~\eqref{eq:Pprimaldual} has a solution, which is a saddle-point of~$\mathcal{L}$, denoted by~$(\xi^*,y^*)$.
Choose~$\tau>0$,~$\sigma>0$ and~$0<\theta\leq 1$ such that
\begin{equation}
\max\left\{\frac{1}{\tau\tilde\gamma+1},\frac{1}{\sigma\tilde\delta+1}\right\}\leq\theta\leq\frac{1}{L_K^2\tau\sigma}.
\label{eq:conditionstep}
\end{equation}
Then, for any~$\omega$ such that
\begin{equation}
\max\left\{\frac{1}{\tau\tilde\gamma+1},\frac{\theta+1}{\sigma\tilde\delta+2}\right\}\leq\omega\leq\theta
\label{eq:conditionomega}
\end{equation}
we have the following majoration for any~$N\in\mathbb{N}$ and any~$(\xi,y)\in Z\times Y$:
\begin{equation}
\begin{split}
\frac{1}{2\tau}\,\lVert \xi_N - \xi \rVert^2
+(1-\omega L_K^2\tau\sigma)\frac{1}{2\sigma}\,\lVert y_N - y\rVert^2
+\sum_{n=1}^N \frac{\omega^n}{\omega^{n-1}}\big(\mathcal{L}(\xi_n;y) - \mathcal{L}(\xi;y_n)\big)
\\[2mm]\leq
\frac{\omega^N}{2\tau}\,\lVert \xi_0 - \xi \rVert^2 + \frac{\omega^N}{2\sigma}\,\lVert y_0 - y \rVert^2
\end{split}
\label{eq:CVgap0}
\end{equation}
where $(\xi_n,y_n)_n$ are generated by Algorithm~\eqref{eq:algoADMM-PD22}.
Hence, if we define 
\begin{equation}
T_N 
:= \sum_{n=1}^N \frac{1}{\omega^{n-1}}
=\frac{1-\omega^N}{\omega^{N-1}(1-\omega)}
\label{eq:defTN}
\end{equation}
and let
\begin{equation}
\Xi_N := \frac{1}{T_N}\sum_{n=1}^N\frac{1}{\omega^{n-1}} \, \xi_n
\qquad\text{and}\qquad
Y_N := \frac{1}{T_N}\sum_{n=1}^N\frac{1}{\omega^{n-1}} \, y_n.
\end{equation}
Then we have the following bound for any~$(\xi,y)\in Z\times Y$:
\begin{equation}
\begin{split}
\frac{1-\omega}{\omega(1-\omega^N)}\,\frac{1}{2\tau}\,\lVert \xi-\xi_N \rVert^2&
+\frac{1-\omega}{\omega(1-\omega^N)}\,(1-\omega L_K^2\tau\sigma)\frac{1}{2\sigma}\,\lVert y-y_N\rVert^2\\[2mm]
&+\mathcal{L}(\Xi_N;y) - \mathcal{L}(\xi;Y_N)\\[2mm]&
\leq
\frac{1}{T_N}\,\frac{1}{2\tau}\,\lVert \xi-\xi_0 \rVert^2
 + \frac{1}{T_N}\,\frac{1}{2\sigma}\,\lVert y-y_0 \rVert^2.
\end{split}
\label{eq:bound}
\end{equation}
\end{theorem}

This theorem provides a linear ergodic convergence rate, namely for the sequences~$(\Xi_N)$ and~$(Y_N)$. This rate can be compared with~\cite{nesterov2004introductory}, and will proved to be better with optimal parameters. Also note that no assumption is made about the rank of the linear operator~$K$. 
Equation~\eqref{eq:bound} can be applied to~$\xi=\xi^*$ and~$y=y^*$, which yields a nonergodic linear convergence rate for the variable convergence (see subsection~\ref{sub:linear}).

A similar result may be found in~\cite{chambolle2011first}, but the rate we provide here is better, since no restrictive assumptions are made on the parameters values, unless necessary.

\subsection{Proof of convergence}
We proceed analogously to the proof in~\cite{chambolle2011first}, but we do not specify any parameter unless needed. This proof is also inspired by the one found in \cite{chambolle2015ergodic}, which does not allow $\theta\neq 1$. For now, we only assume that~$0<\theta\leq1$.
\medskip

\subsubsection{Preliminaries}
Let us consider the general updates $(\hat y,\hat \xi)$ by setting for any~$(\bar \xi,\tilde \xi)\in Z^2$ and~$(\bar y,\tilde y)\in Y^2$
\begin{equation}
\begin{cases}
\hat y := \text{prox}_{\sigma H^*}(\bar y + \sigma K\tilde \xi)\\[2mm]
\hat \xi := \text{prox}_{\tau G}(\bar \xi - \tau K^*\tilde y).
\end{cases}
\end{equation}
In other terms, $\hat y$ and $\hat\xi$ are the output of an iteration, and are respectively computed from initial points $(\bar y,\tilde\xi)$ and $(\tilde y,\bar\xi)$. 
These points are related by first-order optimality conditions.
For instance, the point~$\hat \xi$ is defined as the solution of a minimization problem
\begin{equation}
\hat \xi = \arg\min_{\xi\in Z} \left\{ \frac{1}{2\tau}\,\lVert\bar \xi - \tau K^*\tilde y -\xi\rVert^2 + G(\xi)\right\}
\end{equation}
so, by optimality, we obtain
\begin{equation}
-\frac{1}{\tau}\,(\hat \xi-\bar \xi)  - K^*\tilde y \in \partial G(\hat \xi).
\end{equation}
Similarly, the definition of~$\hat y$ yields
\begin{equation}
-\frac{1}{\sigma}\,(\hat y-\bar y)  + K\tilde \xi \in \partial F^*(\hat y).
\end{equation}

Using the definition of strong convexity recalled in~\eqref{eq:fortconvex}, we get (after expanding the scalar products)
\begin{equation}
G(\xi) + \frac{1}{2\tau}\,\lVert \xi-\bar \xi\rVert^2
 \geq G(\hat \xi) + \langle K(\hat \xi-\xi),\tilde y\rangle + \frac{1}{2\tau}\,\lVert \hat \xi-\bar \xi\rVert^2 + \frac{1}{2\tau}\,\lVert \xi-\hat \xi\rVert^2 + \frac{\tilde\gamma}{2}\,\lVert \xi-\hat \xi\rVert^2
\label{eq:Ggamma}
\end{equation}
\begin{equation}
H^*(y) + \frac{1}{2\sigma}\,\lVert y-\bar y\rVert^2
 \geq H^*(\hat y) - \langle K\tilde \xi,\hat y - y\rangle + \frac{1}{2\sigma}\,\lVert \hat y-\bar y\rVert^2 + \frac{1}{2\sigma}\,\lVert y-\hat y\rVert^2 + \frac{\tilde\delta}{2}\,\lVert y-\hat y\rVert^2.
\label{eq:Fdelta}
\end{equation} 
Now, summing~\eqref{eq:Ggamma} and~\eqref{eq:Fdelta}, we have after rearrangement
\begin{equation}
\begin{split}
\mathcal{L}(\hat \xi;y) - \mathcal{L}(\xi;\hat y)
\leq&\ 
\frac{1}{2\tau}\,\lVert \xi-\bar \xi\rVert^2
 - \frac{1+\tau\tilde\gamma}{2\tau}\,\lVert \xi-\hat \xi\rVert^2
 - \frac{1}{2\tau}\,\lVert\bar \xi-\hat \xi\rVert^2\\
& + \frac{1}{2\sigma}\,\lVert y-\bar y\rVert^2
 - \frac{1+\sigma\tilde\delta}{2\sigma}\,\lVert y-\hat y\rVert^2
 - \frac{1}{2\sigma}\,\lVert\bar y-\hat y\rVert^2\\[2mm]
&+\langle K(\hat \xi-\xi),\hat y -\tilde y\rangle
-\langle K(\hat \xi-\tilde \xi),\hat y - y\rangle.
\end{split} 
\label{somme0}
\end{equation}

\subsubsection{First inequality}
Let us now prove the following lemma:
\begin{lemma}\label{lemma1}
Let $(\xi_n,y_n)_n$ be generated by Algorithm~\eqref{eq:algoADMM-PD22}. Then, for any $n\in\mathbb{N}$, $\tau,\sigma>0$ and $0<\omega\leq\theta$, we have
\begin{equation}
\begin{split}
\mathcal{L}(\xi_n;y) - \mathcal{L}(\xi;y_n)
\leq&\
 \frac{1}{2\tau}\,\lVert \xi-\xi_n\rVert^2 + \frac{1}{2\sigma}\,\lVert y-y_n\rVert^2\\[2mm]
& -\frac{1}{\omega}\left( \frac{1}{2\tau}\,\lVert \xi-\xi_{n+1}\rVert^2 + \frac{1}{2\sigma}\,\lVert y-y_{n+1}\rVert^2\right)\\[2mm]
&+\omega\,\frac{1}{2\tau}\,\lVert \xi_{n-1}-\xi_n\rVert^2
 - \frac{1}{2\tau}\,\lVert \xi_n-\xi_{n+1}\rVert^2\\[2mm]
&+\omega\,\langle K(\xi_{n-1}-\xi_n),y-y_n\rangle
 -\langle K(\xi_n-\xi_{n+1}),y-y_{n+1}\rangle.
\end{split} 
\label{eq:majgapomega0}
\end{equation}
\end{lemma}

\paragraph{Proof}
We specify the six variables in~\eqref{somme0}, by choosing on one hand
\begin{equation}
\hat \xi = \xi_{n+1},\quad
\bar \xi = \xi_n,\quad\text{and}\quad
\tilde \xi = \xi_n + \theta\,(\xi_n-\xi_{n-1})
\end{equation}
for~$1\geq\theta >0$ not specified yet, and~
\begin{equation}
\hat y = y_{n+1},\quad
\bar y = y_n,\quad\text{and}\quad
\tilde y = y_{n+1}
\end{equation}
on the other hand, which leads to the iterations in~\eqref{eq:algoADMM-PD22}. After a simplification, we get
\begin{equation}
\begin{split}
\mathcal{L}(\xi_{n+1};y) - \mathcal{L}(\xi;y_{n+1})
\leq&\
\frac{1}{2\tau}\,\lVert \xi-\xi_n\rVert^2
 + \frac{1}{2\sigma}\,\lVert y-y_n\rVert^2\\
& 
 - \frac{{1+\tau\tilde\gamma}}{2\tau}\,\lVert \xi-\xi_{n+1}\rVert^2
 - \frac{{1+\sigma\tilde\delta}}{2\sigma}\,\lVert y-y_{n+1}\rVert^2\\
& - \frac{1}{2\tau}\,\lVert \xi_n-\xi_{n+1}\rVert^2
 - \frac{1}{2\sigma}\,\lVert y_n-y_{n+1}\rVert^2\\[2mm]
&+\theta\,\langle K(\xi_{n-1}-\xi_n),y-y_{n+1}\rangle\\[2mm]
&
-\langle K(\xi_n-\xi_{n+1}),y-y_{n+1}\rangle.
\label{eq:somme}
\end{split} 
\end{equation}
Now, we define~$\tau\tilde\gamma = \mu>0$ and~$\sigma\tilde\delta = \mu'>0$. For any~$n\in\mathbb{N}$, we set
\begin{equation}
\Delta_n = \frac{1}{2\tau}\,\lVert \xi-\xi_n\rVert^2 + \frac{1}{2\sigma}\,\lVert y-y_n\rVert^2.
\end{equation}
Hence, we can rewrite~\eqref{eq:somme} with~$\Delta_n$, which yields
\begin{equation}
\begin{split}
\mathcal{L}(\xi_{n+1};y) - \mathcal{L}(\xi;y_{n+1})
\leq&\
\Delta_n
 - (1+\mu)\,\Delta_{n+1}\\[2mm]
& - \frac{1}{2\tau}\,\lVert \xi_n-\xi_{n+1}\rVert^2
 - \frac{1}{2\sigma}\,\lVert y_n-y_{n+1}\rVert^2\\[2mm]
&+\theta\,\langle K(\xi_{n-1}-\xi_n),y-y_{n+1}\rangle\\[2mm]
&
 -\langle K(\xi_n-\xi_{n+1}),y-y_{n+1}\rangle
\\[2mm]
& +\frac{\mu-\mu'}{2\sigma}\,\lVert y-y_{n+1}\rVert^2.
\end{split} 
\label{eq:gapDelta'}
\end{equation}
Let us bound the scalar products in~\eqref{eq:gapDelta'}. For any~$0<\omega\leq \theta$, we have the decomposition
\begin{equation}
\begin{split}
\theta\,\langle K(\xi_{n-1}-\xi_n),y-y_{n+1}\rangle
= &\ 
\omega\,\langle K(\xi_{n-1}-\xi_n),y-y_n\rangle\\[2mm]&
+\omega\,\langle K(\xi_{n-1}-\xi_n),y_n-y_{n+1}\rangle\\[2mm]&
+(\theta-\omega)\,\langle K(\xi_{n-1}-\xi_n),y-y_{n+1}\rangle.
\end{split}
\end{equation}
Let us have a closer look at the last two terms. Let~$\alpha>0$. Since~$\omega\geq 0$, we have
\begin{equation}
\begin{split}
\omega\,\langle K(\xi_{n-1}-\xi_n),y_n-y_{n+1}\rangle
&\leq \omega\,L_K\,\lVert \xi_{n-1}-\xi_n\rVert\cdot\lVert y_n-y_{n+1}\rVert\\
&\leq \omega\,L_K\,\left(\frac{\alpha}{2}\,\lVert \xi_{n-1}-\xi_n\rVert^2
+\frac{1}{2\alpha}\,\lVert y_n-y_{n+1}\rVert^2\right).
\label{eq:scalar1}
\end{split}
\end{equation}
Similarly, since~$\theta-\omega\geq0$,
\begin{equation}
(\theta-\omega)\,\langle K(\xi_{n-1}-\xi_n),y-y_{n+1}\rangle
\leq (\theta-\omega)\,L_K\,\left(\frac{\alpha}{2}\,\lVert \xi_{n-1}-\xi_n\rVert^2
+\frac{1}{2\alpha}\,\lVert y-y_{n+1}\rVert^2\right).
\label{eq:scalar2}
\end{equation}
After simplification, the majoration~\eqref{eq:gapDelta'} becomes, thanks to inequalities~\eqref{eq:scalar1} and~\eqref{eq:scalar2},
\begin{equation}
\begin{split}
\mathcal{L}(\xi_{n+1};y) - \mathcal{L}(\xi;y_{n+1})
&\leq
\Delta_n
 - (1+\mu)\,\Delta_{n+1}\\[2mm]
&+\theta\,L_K\,\frac{\alpha}{2}\,\lVert \xi_{n-1}-\xi_n\rVert^2
 - \frac{1}{2\tau}\,\lVert \xi_n-\xi_{n+1}\rVert^2\\[1mm]
&
+\left(\frac{\omega\,L_K}{2\alpha}
 - \frac{1}{2\sigma}\right)\,\lVert y_n-y_{n+1}\rVert^2\\[1mm]
&+\omega\,\langle K(\xi_{n-1}-\xi_n),y-y_n\rangle
 -\langle K(\xi_n-\xi_{n+1}),y-y_{n+1}\rangle\\[1mm]&
+\left(\frac{(\theta-\omega)\,L_K}{2\alpha}+\frac{\mu-\mu'}{2\sigma}\right)\lVert y-y_{n+1}\rVert^2.
\end{split} 
\label{K*}
\end{equation}
Choose~$\alpha = \omega L_K\sigma$. Hence, we have~$\omega L_K/\alpha=1/\sigma$, so that the~$\lVert y_n-y_{n+1}\rVert^2$ term cancels. This leads to:
\begin{equation}
\begin{split}
\mathcal{L}(\xi_{n+1};y) - \mathcal{L}(\xi;y_{n+1})
&\leq
\Delta_n
 - (1+\mu)\,\Delta_{n+1}\\[2mm]
&
 + \omega\,\frac{\theta L_K^2\tau\sigma}{2\tau}\,\lVert \xi_{n-1}-\xi_n\rVert^2
 - \frac{1}{2\tau}\,\lVert \xi_n-\xi_{n+1}\rVert^2\\[2mm]
&+\omega\,\langle K(\xi_{n-1}-\xi_n),y-y_n\rangle
 -\langle K(\xi_n-\xi_{n+1}),y-y_{n+1}\rangle\\[2mm]&
+\left(\frac{\theta-\omega}{\omega}+\mu-\mu'\right)\frac{1}{2\sigma}\,\lVert y-y_{n+1}\rVert^2.
\end{split} 
\label{eq:gapPS}
\end{equation}
Since~$1+\mu = 1/\omega + 1+\mu-1/\omega$, we have
\begin{equation}
 - (1+\mu)\,\Delta_{n+1}
 =  -\frac{1}{\omega}\,\Delta_{n+1}
  + \left(\frac{1}{\omega}-\mu-1\right)\,\left(\frac{1}{2\tau}\,\lVert \xi-\xi_{n+1} \rVert^2+\frac{1}{2\sigma}\,\lVert y-y_{n+1} \rVert^2\right)
\end{equation}
so the right-hand side of~\eqref{eq:gapPS} becomes
\begin{equation}
\begin{split}
\Delta_n-\frac{1}{\omega}\,\Delta_{n+1}&+\omega\,\frac{{\theta L_K^2\tau\sigma}}{2\tau}\,\lVert \xi_n-\xi_{n-1}\rVert^2
 - \frac{1}{2\tau}\,\lVert \xi_n-\xi_{n+1}\rVert^2\\[2mm]
&+\omega\,\langle K(\xi_{n-1}-\xi_n),y-y_n\rangle
 -\langle K(\xi_n-\xi_{n+1}),y-y_{n+1}\rangle\\[2mm]
& + {\left(\frac{1}{\omega}-\mu-1\right)}\frac{1}{2\tau}\,\lVert \xi-\xi_{n+1} \rVert^2\\[2mm]
&
+ {\left(\frac{\theta-\omega}{\omega}+\frac{1}{\omega}-\mu'-1\right)}\frac{1}{2\sigma}\,\lVert y-y_{n+1} \rVert^2.
\end{split} 
\label{eq:majgap}
\end{equation}

It is now time to set conditions on~$\omega$,~$\theta$,~$\tau$ and~$\sigma$. First, choose~$\theta$,~$\tau$ and~$\sigma$ so that~$\theta L_K^2\tau\sigma \leq 1$. Then, choose~$\theta$ so that both~$1/\omega-\mu-1$ and~$(\theta-\omega)/\omega+1/\omega-\mu'-1$ are nonpositive, which implies that
\begin{equation}
\frac{1}{\mu+1}\leq \omega \leq \theta
\qquad\text{and}\qquad
\frac{\theta+1}{\mu'+2} \leq \omega \leq \theta.
\end{equation}
Then we can bound~\eqref{eq:majgap} by
\begin{equation}
\begin{split}
\Delta_n-\frac{1}{\omega}\,\Delta_{n+1}
&+\omega\,\frac{1}{2\tau}\,\lVert \xi_{n-1}-\xi_n\rVert^2
 - \frac{1}{2\tau}\,\lVert \xi_n-\xi_{n+1}\rVert^2\\[2mm]
&+\omega\,\langle K(\xi_{n-1}-\xi_n),y-y_n\rangle
 -\langle K(\xi_n-\xi_{n+1}),y-y_{n+1}\rangle.
\end{split} 
\label{eq:majgapomega1}
\end{equation}

\noindent Eventually, back to~\eqref{eq:gapPS} we get the wanted inequality
\begin{equation}
\begin{split}
\mathcal{L}(\xi_{n+1};y) - \mathcal{L}(\xi;y_{n+1})
\leq&\ 
\Delta_n-\frac{1}{\omega}\,\Delta_{n+1}\\[0mm]
&+\omega\,\frac{1}{2\tau}\,\lVert \xi_{n-1}-\xi_n\rVert^2
 - \frac{1}{2\tau}\,\lVert \xi_n-\xi_{n+1}\rVert^2\\[2mm]
&+\omega\,\langle K(\xi_{n-1}-\xi_n),y-y_n\rangle\\[2mm]
&
 -\langle K(\xi_n-\xi_{n+1}),y-y_{n+1}\rangle.~
\end{split} 
\label{eq:majgapomega}
\end{equation}

\subsubsection{Linear convergence of the iterates}\label{sub:linear}
Multiplying~\eqref{eq:majgapomega} by~$1/\omega^n$ and summing between~$n=0$ and~$n=N-1$ (choose~$\xi^{-1} = \xi^0$) cancels most of the terms:
\begin{equation}
\begin{split}
\sum_{n=1}^N\frac{1}{\omega^{n-1}}\big(\mathcal{L}(\xi_n;y) - \mathcal{L}(\xi;y_n)\big)
&\leq
\Delta_0-\frac{1}{\omega^N}\,\Delta_N
 - \frac{1}{2\tau\omega^{N-1}}\,\lVert \xi_{N-1}-\xi_N\rVert^2\\
&\qquad -\frac{1}{\omega^{N-1}}\,\langle K(\xi_{N-1}-\xi_N),y-y_N\rangle.
\label{eq:telescope}
\end{split}
\end{equation}
Once again, we bound the scalar product: let~$\beta>0$,
\begin{equation}
-\frac{1}{\omega^{N-1}}\,\langle K(\xi_{N-1}-\xi_N),y-y_N\rangle
\leq \frac{L_K}{\omega^{N-1}} \left(\frac{\beta}{2}\,\lVert \xi_{N-1}-\xi_N \rVert^2
+\frac{1}{2\beta}\,\lVert y-y_N\rVert^2\right)
\end{equation}
and inequality~\eqref{eq:telescope} becomes
\begin{equation}
\begin{split}
\sum_{n=1}^N\frac{1}{\omega^{n-1}}\big(\mathcal{L}(\xi_n;y) - \mathcal{L}(\xi;y_n)\big)
\leq&\
\Delta_0-\frac{1}{\omega^N}\,\Delta_N\\&
+ {\left(\frac{L_K\beta}{2\omega^{N-1}} - \frac{1}{2\tau\omega^{N-1}}\right)}\,\lVert \xi_{N-1}-\xi_N\rVert^2\\&+ \frac{L_K}{\omega^{N-1}}\frac{1}{2\beta}\,\lVert y_N-y\rVert^2.
 \end{split}
\end{equation}
Now choose~$\beta = 1/(L_K\tau)$, which cancels the~$\lVert \xi_{N-1}-\xi_N\rVert^2$ term, and we get
\begin{equation}
\sum_{n=1}^N\frac{1}{\omega^{n-1}}\big(\mathcal{L}(\xi_n;y) - \mathcal{L}(\xi;y_n)\big)
\leq
\Delta_0-\frac{1}{\omega^N}\,\Delta_N
 + \frac{L_K^2\tau\sigma}{\omega^{N-1}}\frac{1}{2\sigma}\,\lVert y-y_N\rVert^2.
\label{eq:majgapDelta'}
\end{equation}
Replacing~$\Delta_0$ and~$\Delta_n$ by their respective definition, we obtain
\begin{equation}
\begin{split}
\sum_{n=1}^N\frac{1}{\omega^{n-1}}\big(\mathcal{L}(\xi_n;y) - \mathcal{L}(\xi;y_n)\big)
\leq&\
\frac{1}{2\tau}\,\lVert \xi-\xi_0 \rVert^2 + \frac{1}{2\sigma}\,\lVert y-y_0 \rVert^2\\
&-\frac{1}{\omega^N}\,\frac{1}{2\tau}\,\lVert \xi-\xi_N \rVert^2\\[2mm]
&
-\frac{1}{\omega^N}\,{(1-\omega L_K^2\tau\sigma)}\frac{1}{2\sigma}\,\lVert y-y_N\rVert^2.
\end{split}
\end{equation}
Since~$\omega L_K^2\tau\sigma \leq \theta L_K^2\tau\sigma \leq 1$ and~$\mathcal{L}(\xi_n;y) - \mathcal{L}(\xi;y_n)\geq 0$ for any~$n\in\mathbb{N}$, we have
\begin{equation}
\begin{split}
0\leq \frac{1}{2\tau}\,\lVert \xi-\xi_N \rVert^2
+(1-\omega L_K^2\tau\sigma)\frac{1}{2\sigma}\,\lVert y-y_N\rVert^2
 +\sum_{n=1}^N\frac{\omega^N}{\omega^{n-1}}\big(\mathcal{L}(\xi_n;y) - \mathcal{L}(\xi;y_n)\big)
\\[2mm]\leq
\frac{\omega^N}{2\tau}\,\lVert \xi-\xi_0 \rVert^2 + \frac{\omega^N}{2\sigma}\,\lVert y-y_0 \rVert^2.
\end{split}
\label{eq:CVgap}
\end{equation}

The latter inequality proves the linear convergence of the iterates:
\begin{corollary}\label{lemma2}
Assume problem~\eqref{eq:Pprimaldual} has a solution, which is a saddle-point of~$\mathcal{L}$, denoted by~$(\xi^*,y^*)$.
Let $(\xi_n,y_n)_n$ be generated by Algorithm~\eqref{eq:algoADMM-PD22}. Suppose there exist $\tau$, $\sigma$, $\theta$ and $\omega$ satisfying both conditions~\eqref{eq:conditionstep} and~\eqref{eq:conditionomega}. Then, for any $N\in\mathbb{N}$, we have
\begin{equation}
\lVert \xi^*-\xi_N \rVert^2
\leq
\omega^N\left(\lVert \xi^*-\xi_0 \rVert^2 + \frac{\tau}{\sigma}\,\lVert y^*-y_0 \rVert^2\right).
\end{equation}
Moreover, if $\omega L_K^2\tau\sigma\neq 1$, then we also have
\begin{equation}
\lVert y^*-y_N\rVert^2
\leq
\frac{\omega^N}{1-\omega L_K^2\tau\sigma}\left(\frac{\sigma}{\tau}\lVert \xi^*-\xi_0 \rVert^2 + \lVert y^*-y_0 \rVert^2\right).
\end{equation}
\end{corollary}

\noindent\textsc{Remark:} The convergence rate in Corollary~\ref{lemma2} can be improved if we use the fact that, by definition,
\begin{equation}\begin{split}
\mathcal{L}(\xi_{n+1};y^*) - \mathcal{L}(\xi^*;y_{n+1})
 = G(\xi_{n+1}) - G(\xi^*)
 ~&+ H^*(y_{n+1}) - H^*(y^*)\\
 &+ \langle K\xi_{n+1},y^*\rangle - \langle K\xi^*,y_{n+1}\rangle.
 \end{split}
\end{equation}
The strong convexity of $G$ et $H^*$ and the optimality of $\xi^*$ and $y^*$ yield the following inequalities:
\begin{eqnarray}
G(\xi_{n+1}) - G(\xi^*) \geq \langle -K^*y^*,\xi_{n+1}-\xi^*\rangle + \frac{\tilde\gamma}{2}\,\lVert\xi_{n+1}-\xi^*\rVert^2\\
H^*(y_{n+1}) - H^*(y^*) \geq \langle K\xi^*,y_{n+1}-y^*\rangle + \frac{\tilde\delta}{2}\,\lVert y_{n+1}-y^*\rVert^2.
\end{eqnarray}
This implies that
\begin{equation}
      \frac{\tilde\gamma}{2}\,\lVert\xi_{n+1}-\xi^*\rVert^2
     + \frac{\tilde\delta}{2}\,\lVert y_{n+1}-y^*\rVert^2
    \leq 
\mathcal{L}(\xi_{n+1};y^*) - \mathcal{L}(\xi^*;y_{n+1})
\end{equation}
since the sum of the scalar products cancels. Hence, if we choose not to control the primal-dual gap, choosing $(\xi,y) = (\xi^*,y^*)$, in~\eqref{eq:somme} becomes
\begin{equation}
\begin{split}
0
\leq\ &
\frac{1}{2\tau}\,\lVert \xi^*-\xi_n\rVert^2
 + \frac{1}{2\sigma}\,\lVert y^*-y_n\rVert^2\\
& 
 - \frac{{1+2\tau\tilde\gamma}}{2\tau}\,\lVert \xi^*-\xi_{n+1}\rVert^2
 - \frac{{1+2\sigma\tilde\delta}}{2\sigma}\,\lVert y^*-y_{n+1}\rVert^2\\
& - \frac{1}{2\tau}\,\lVert \xi_n-\xi_{n+1}\rVert^2
 - \frac{1}{2\sigma}\,\lVert y_n-y_{n+1}\rVert^2\\[2mm]
&+\theta\,\langle K(\xi_{n-1}-\xi_n),y^*-y_{n+1}\rangle
-\langle K(\xi_n-\xi_{n+1}),y^*-y_{n+1}\rangle
\label{eq:somme2}
\end{split} 
\end{equation}
which means that all the computations from~\eqref{eq:somme} to~\eqref{eq:CVgap} hold, with $\mu$ and $\mu'$ replaced by $\tilde\mu = 2\mu$ and $\tilde\mu' = 2\mu'$ and without $\mathcal{L}$-terms, as well as the constraints on the parameters. In others terms, the same computations prove that
\begin{corollary}\label{lemma2new}
Assume problem~\eqref{eq:Pprimaldual} has a solution, which is a saddle-point of~$\mathcal{L}$, denoted by~$(\xi^*,y^*)$.
Let $(\xi_n,y_n)_n$ be generated by Algorithm~\eqref{eq:algoADMM-PD22}. Suppose there exist $\tau$, $\sigma$, $\theta$ and $\omega$ satisfying both conditions
\begin{equation}
\max\left\{\frac{1}{2\tau\tilde\gamma+1},\frac{1}{2\sigma\tilde\delta+1}\right\}\leq\theta\leq\frac{1}{L_K^2\tau\sigma}.
\label{eq:conditionstep2}
\end{equation}
Then, for any~$\tilde\omega$ such that
\begin{equation}
\max\left\{\frac{1}{2\tau\tilde\gamma+1},\frac{\theta+1}{2\sigma\tilde\delta+2}\right\}\leq\tilde\omega\leq\theta.
\label{eq:conditionomega2}
\end{equation}
Then, for any $N\in\mathbb{N}$, we have
\begin{equation}
\lVert \xi^*-\xi_N \rVert^2
\leq
\tilde\omega^N\left(\lVert \xi^*-\xi_0 \rVert^2 + \frac{\tau}{\sigma}\,\lVert y^*-y_0 \rVert^2\right).\label{eq:CVgap2}
\end{equation}
Moreover, if $\tilde\omega L_K^2\tau\sigma\neq 1$, then we also have
\begin{equation}
\lVert y^*-y_N\rVert^2
\leq
\frac{\tilde\omega^N}{1-\tilde\omega L_K^2\tau\sigma}\left(\frac{\sigma}{\tau}\lVert \xi^*-\xi_0 \rVert^2 + \lVert y^*-y_0 \rVert^2\right).
\end{equation}
\end{corollary}
For given $\tau$, $\sigma$ and $\theta$, the lower bounds $1/(2\tau\tilde\gamma+1)$ and $(\theta+1)/(2\sigma\tilde\delta+2)$ for $\tilde\omega$ are smaller than those for $\omega$. Thus, the new rate $\tilde\omega$ can be expected to be better than the global one $\omega$ (which is called \emph{global} since it also holds for the objectif error, as shown in the next paragraph). This will be checked in Subsection~\ref{sub:paramPDHG}.

\subsubsection{End of the proof}
We can now complete the proof of Theorem~\ref{thm1}.
Dividing~\eqref{eq:CVgap} by~$\omega^N\neq0$ and by~$T_N\neq0$, we get
\begin{equation}
\begin{split}
\frac{1-\omega}{\omega(1-\omega^N)}\,\frac{1}{2\tau}\,\lVert \xi-\xi_N \rVert^2
&+\frac{1-\omega}{\omega(1-\omega^N)}\,(1-\omega L_K^2\tau\sigma)\frac{1}{2\sigma}\,\lVert y-y_N\rVert^2\\
&
+\frac{1}{T_N}\,\sum_{n=1}^N\frac{1}{\omega^{n-1}}\big(\mathcal{L}(\xi_n;y) - \mathcal{L}(\xi;y_n)\big)\\
&\leq
\frac{1}{T_N}\,\frac{1}{2\tau}\,\lVert \xi-\xi_0 \rVert^2 + \frac{1}{T_N}\,\frac{1}{2\sigma}\,\lVert y-y_0 \rVert^2.
\label{eq:ergo}
\end{split}
\end{equation}
But, by convexity,
\begin{equation}
\mathcal{L}(\Xi_N;y) - \mathcal{L}(\xi;Y_N)
\leq
\frac{1}{T_N}\,\sum_{n=1}^N\frac{1}{\omega^{n-1}}\big(\mathcal{L}(\xi_n;y) - \mathcal{L}(\xi;y_n)\big)
\end{equation}
Therefore,~\eqref{eq:ergo} becomes
\begin{equation}
\begin{split}
 \frac{1-\omega}{\omega(1-\omega^N)}\,\frac{1}{2\tau}\,\lVert \xi-\xi_N \rVert^2&
+\frac{1-\omega}{\omega(1-\omega^N)}\,(1-\omega L_K^2\tau\sigma)\frac{1}{2\sigma}\,\lVert y-y_N\rVert^2\\[2mm]&
+\mathcal{L}(\Xi_N;y) - \mathcal{L}(\xi;Y_N)\\[2mm]&
\leq
\frac{1}{T_N}\,\frac{1}{2\tau}\,\lVert \xi-\xi_0 \rVert^2
 + \frac{1}{T_N}\,\frac{1}{2\sigma}\,\lVert y-y_0 \rVert^2
\end{split}
\label{eq:ergoobj}
\end{equation}
which completes the proof of Theorem~\ref{thm1}.

\paragraph{Remark}
Equation \eqref{eq:ergoobj} provides a way to establish the (ergodic) linear convergence of the objective function whe applied to $(\xi,y)=(\xi^*_N,y^*_N)$ where the supremum of the primal-dual gap is attained. The main argument relies on the Lipschitz continuity of the gradients of $H$ and $G^*$, which ensures that this point is close to $(\xi^*,y^*)$. A more detailed example of such computations is provided in Section \ref{sec:cvADMM} in the case of the ADMM.

\subsection{Choice of parameters}\label{sub:paramPDHG}
Theorem~\ref{thm1} holds provided one can properly choose the steps~$\tau$ and~$\sigma$ and the relaxation parameter~$\theta$. We study some particular choices for those parameters and the convergence rate they yield. Since a smaller~$\omega$ leads to a faster convergence, we tune the algorithm parameters to minimize the lower bound of~$\omega$. Here is how we proceed:

\begin{enumerate}
\item Fix~$\tau>0$.
\item Find conditions on~$\sigma$ so that inequalities~\eqref{eq:conditionstep} hold.
\item Minimize~$(\theta+1)/(\sigma\tilde\delta+2)$ with respect to (w.r.t.)~$\theta$ satisfying~\eqref{eq:conditionstep} and w.r.t~$\sigma$ given by the previous step.
\item Compare this minimum to~$1/(\tau\tilde\gamma+1)$ and deduce the lower bound~$\omega^*(\tau)$ for~$\omega$.
\item Minimize~$\omega^*(\tau)$ and derive the optimal rate~$\omega^*$.
\end{enumerate}
Since the resulting parameters are compatible with conditions~\eqref{eq:conditionstep2} and~\eqref{eq:conditionomega2}, the left-hand member in~\eqref{eq:CVgap2} yields a better theoretical rate for the convergence of the variables.

Besides, the same computations (with $\tilde\gamma$ and $\tilde\delta$ doubled) may be used to choose the parameters so that the rate $\tilde\omega$ (Corollary~\ref{lemma2new}) is minimal.

\subsubsection{Case~$\theta=1$}
We first fix~$\theta=1$. As shown in~\cite{chambolle2011first}, this choice is equivalent to the ADMM with an additional proximal term.

Fix~$\tau>0$. Replacing~$\theta=1$ in~\eqref{eq:conditionstep}, we obtain that the steps~$\tau$ and~$\sigma$ are constrained as following
\begin{equation}
1\leq \frac{1}{L_K^2\tau\sigma}
\label{eq:theta1sigma}
\end{equation}
which implies that~$\sigma\leq 1/(L_K^2\tau)$.
Then,~\eqref{eq:conditionomega} in Theorem 1 states that the convergence rate~$\omega$ satisfies
\begin{equation}
\max\left\{\frac{1}{\tau\tilde\gamma+1},\frac{1}{\sigma\tilde\delta/2+1}\right\}
\leq \omega\leq 1
\end{equation}
Let us minimize~$1/(\sigma\tilde\delta/2+1)$ w.r.t.~$\sigma$ satisfying~\eqref{eq:theta1sigma}. Since the map~$\sigma\mapsto1/(\sigma\tilde\delta/2+1)$ is nondecreasing, its minimum is reached when~$\sigma$ is maximal, which leads to 
\begin{equation}
\min_{\sigma\text{ subject to }\eqref{eq:theta1sigma}}\left\{\frac{1}{\sigma\tilde\delta/2+1}\right\} = 
\frac{1}{\tilde\delta/(2L_K^2\tau)+1}.
\end{equation}
Now, compare it to~$1/(\tau\tilde\gamma+1)$. It is clear that the quantity~$1/(\tilde\delta/(2L_K^2\tau)+1)$ is greater than~$1/(\tau\tilde\gamma+1)$ as soon as~$\tau^2\geq\tilde\delta/(2\tilde\gamma L_K^2)$. Hence, the lower bound~$\omega^*(\tau)$ is given by
\begin{equation}
\omega^*(\tau)
= \max\left\{\frac{1}{\tau\tilde\gamma+1},\frac{1}{\tilde\delta/(2L_K^2\tau)+1}\right\}
= \begin{cases}
\displaystyle \frac{1}{\tau\tilde\gamma+1} & \text{if }0<\tau<\sqrt{\tilde\delta/(2\tilde\gamma L_K^2)}\\[3mm]
\displaystyle \frac{1}{\tilde\delta/(2L_K^2\tau)+1} & \text{if }\tau\geq\sqrt{\tilde\delta/(2\tilde\gamma L_K^2)}
\end{cases}
\end{equation}
which is minimal for~$\tau^*=\sqrt{\tilde\delta/(2\tilde\gamma L_K^2)}$ and leads to the optimal rate
\begin{equation}
\omega^*
 = \omega^*(\tau^*)
 = \frac{1}{\sqrt{(\tilde\gamma\tilde\delta)/(2L_K^2)}+1}
 = \frac{1}{\sqrt{1/(2\kappa_F)}+1}.
\end{equation}
This rate is reached for
\begin{equation}
\tau = \tau^* = \sqrt{\frac{\tilde\delta}{2\tilde\gamma L_K^2}}
\quad\text{and}\quad
\sigma  = \frac{1}{L_K^2\tau^*} = \sqrt{\frac{2\tilde\gamma}{\tilde\delta L_K^2}}.
\end{equation}
One can check that the same choice for $\tau$ and $\sigma$ yield the minimal value for the solution error rate $\tilde\omega$, which is
\begin{equation}
\tilde\omega^*
 = \frac{1}{2\sqrt{(\tilde\gamma\tilde\delta)/(2L_K^2)}+1}
 = \frac{1}{\sqrt{2/\kappa_F}+1}.
\end{equation}
In other terms, in the case where $\theta=1$, the best choice for the global rate $\omega$ and for the solution error rate $\tilde\omega$ coincide.

\subsubsection{The best convergence rate ($\theta<1$)\label{sub:best}}
In this section, we want to derive the best convergence rate given the constraints in~\eqref{thm1}.
\begin{theorem}\label{thm2}
The best convergence rate in Theorem 1 is obtained when choosing 
\begin{equation}
\tau = \frac{\tilde\delta}{2\,L_K^2}\left(1+\sqrt{1+\frac{4\,L_K^2}{\tilde\gamma\tilde\delta}}\right)
\qquad\text{and}\qquad
\sigma = \frac{\tilde\gamma}{2\,L_K^2}\left(1+\sqrt{1+\frac{4\,L_K^2}{\tilde\gamma\tilde\delta}}\right)
\end{equation}
and, if $\kappa_F = L_K^2/(\tilde\gamma\tilde\delta)$,
\begin{equation}
\theta = 
\frac{\displaystyle \sqrt{1+(4\,L_K^2)/(\tilde\gamma\tilde\delta)}-1}{\displaystyle \sqrt{1+(4\,L_K^2)/(\tilde\gamma\tilde\delta)}+1} 
= \frac{\displaystyle \sqrt{1+4\,\kappa_F}-1}{\displaystyle \sqrt{1+4\,\kappa_F}+1} <1
\end{equation}
which satisfy~$\tau\tilde\gamma=\sigma\tilde\delta$. The resulting rate is $\omega^* = \theta$.
\end{theorem}
\paragraph{Proof} 
Fix~$\tau>0$ and find out which conditions~$\sigma$ must satisfy to ensure the existence of~$\theta$ satisfying ~\eqref{eq:conditionstep}. There exists~$\theta$ satisfying~\eqref{eq:conditionstep} if
\begin{equation}
\frac{1}{\tau\tilde\gamma+1} \leq \frac{1}{L_K^2\tau\sigma}
\quad\text{and}\quad
\frac{1}{\sigma\tilde\delta+1} \leq \frac{1}{L_K^2\tau\sigma}.
\label{eq:conditionssigma}
\end{equation}
which also reads
\begin{equation}
\sigma \leq \frac{1}{L_K^2\tau}+\frac{\tilde\gamma}{L_K^2}
\quad\text{and}\quad
(L_K^2\tau-\tilde\delta)\,\sigma \leq 1.
\end{equation}

Let us determine conditions on~$\sigma$ so that these inequalities hold. If~$L_K^2\tau-\tilde\delta\leq0$, i.e.~$\tau\leq\tilde\delta/L_K^2$, the second inequality is always true. Hence, let us study the case~$L_K^2\tau-\tilde\delta>0$, i.e.~$\tau>\tilde\delta/L_K^2$. It implies that~$\sigma$ must satisfy both majorations
\begin{equation}
\sigma \leq \frac{1}{L_K^2\tau}+\frac{\tilde\gamma}{L_K^2}
\quad\text{and}\quad
\sigma \leq \frac{1}{L_K^2\tau-\tilde\delta}.
\label{eq:condition1}
\end{equation}
Let us compare these two bounds. Since
\begin{equation}
\frac{1}{L_K^2\tau}+\frac{\tilde\gamma}{L_K^2} - \frac{1}{L_K^2\tau-\tilde\delta}
=
\frac{\tilde\gamma L_K^2\tau^2-\tilde\gamma\tilde\delta\tau-\tilde\delta}{L_K^2\tau(L_K^2\tau-\tilde\delta)}
\end{equation}
with~$L_K^2\tau(L_K^2\tau-\tilde\delta)$ positive,~$1/(L_K^2\tau)+\tilde\gamma/L_K^2$ is greater than~$1/(L_K^2\tau-\tilde\delta)$ iff~$\tilde\gamma L_K^2\tau^2-\tilde\gamma\tilde\delta\tau-\tilde\delta \geq 0$, i.e. iff~$\tau\geq\tau^*$, given by
\begin{equation}
\tau^*= \frac{\tilde\delta}{2\,L_K^2}\left(1+\sqrt{1+\frac{4\,L_K^2}{\tilde\gamma\tilde\delta}}\right)
>\frac{\tilde\delta}{L_K^2}.
\label{eq:taustar}
\end{equation}
Therefore, for any~$\tilde\delta/L_K^2<\tau\leq\tau^*$,~\eqref{eq:condition1} becomes~$\sigma \leq 1/(L_K^2\tau)+\tilde\gamma/L_K^2$.
If~$\tau>\tau^*$,~\eqref{eq:condition1} reads~$\sigma \leq 1/(L_K^2\tau-\tilde\delta)$. As a conclusion, we have the following upper bounds for~$\sigma$:
\begin{equation}
\sigma\leq\begin{cases}
\displaystyle\frac{1}{L_K^2\tau}+\frac{\tilde\gamma}{L_K^2} & \text{if } 0<\tau\leq\tau^*\\[2mm]
\displaystyle\frac{1}{L_K^2\tau-\tilde\delta} &\text{if } \tau^*<\tau.
\end{cases}
\label{eq:conditionsigma}
\end{equation}

Now, fix~$\sigma$ satisfying~\eqref{eq:conditionsigma} and let us minimize~$(\theta+1)/(\sigma\tilde\delta+2)$ subject to~\eqref{eq:conditionstep}. The map~$\theta\mapsto(\theta+1)/(\sigma\tilde\delta+2)$ is minimal when~$\theta$ is minimal. Hence, let us determine the lower bound of~$\theta$, which is given by
\begin{equation}
\max\left\{\frac{1}{\tau\tilde\gamma+1},\frac{1}{\sigma\tilde\delta+1}\right\}.
\end{equation}
First, remark that, if~$\tau>\tilde\delta/L_K^2$, then
\begin{equation}
\frac{\tilde\delta}{L_K^2\tau-\tilde\delta} \leq \tau\tilde\gamma
\quad\Longleftrightarrow\quad
\tilde\gamma L_K^2\tau^2-\tilde\gamma\tilde\delta\tau - \tilde\delta \geq 0
\quad\Longleftrightarrow\quad
\tau\geq\tau^*.
\label{eq:remark1}
\end{equation}
Suppose that~$\tau>\tau^*$, which implies that~$\tau>\tilde\delta/L_K^2$. Since~$\sigma$ is bounded from above by~$1/(L_K^2\tau-\tilde\delta)$, we deduce that~$\sigma\tilde\delta\leq\tau\tilde\gamma$, which yields
\begin{equation}
\max\left\{\frac{1}{\tau\tilde\gamma+1},\frac{1}{\sigma\tilde\delta+1}\right\} = 
\frac{1}{\sigma\tilde\delta+1}\quad\text{if}\quad0<\sigma\leq\frac{1}{L_K^2\tau-\tilde\delta}.
\label{eq:thetasigma1}
\end{equation}
Now, let us consider the case~$\tau\leq\tau^*$. Since
\begin{equation}
\frac{1}{L_K^2\tau}+\frac{\tilde\gamma}{L_K^2} \geq \frac{\tau\tilde\gamma}{\tilde\delta}
\quad\Longleftrightarrow\quad
\tilde\gamma L_K^2\tau^2 - \tilde\gamma\tilde\delta\tau - \tilde\delta \leq 0
\quad\Longleftrightarrow\quad
\tau\leq\tau^*
\end{equation}
we deduce that
\begin{equation}
\max\left\{\frac{1}{\tau\tilde\gamma+1},\frac{1}{\sigma\tilde\delta+1}\right\}
=\begin{cases}
\displaystyle \frac{1}{\sigma\tilde\delta+1} & \displaystyle \quad\text{if}\quad 0<\sigma\leq\frac{\tau\tilde\gamma}{\tilde\delta}\\[3mm]
\displaystyle \frac{1}{\tau\tilde\gamma+1} & \displaystyle \quad\text{if}\quad \frac{\tau\tilde\gamma}{\tilde\delta}<\sigma\leq\frac{1}{L_K^2\tau}+\frac{\tilde\gamma}{L_K^2}.
\end{cases}
\label{eq:thetasigma2}
\end{equation}

Let us minimize~$(\theta+1)/(\sigma\tilde\delta+2)$ w.r.t. to~$\sigma$, when~$\theta$ is equal to its lower bound~$\theta^*(\sigma)$, given by~\eqref{eq:thetasigma1} and~\eqref{eq:thetasigma2}. This leads to minimize the following quantity w.r.t.~$\sigma$:
\begin{equation}
\frac{\theta^*(\sigma)+1}{\sigma\tilde\delta+2} =
\begin{cases}
\displaystyle \frac{1}{\sigma\tilde\delta+1} &\displaystyle  \quad\text{if}\quad \tau>\tau^* \text{ or } \left(\tau\leq\tau^* \text{ and } 0<\sigma\leq\frac{\tau\tilde\gamma}{\tilde\delta}\right)\\[3mm]
\displaystyle \frac{1}{\tau\tilde\gamma+1}\,\frac{\tau\tilde\gamma+2}{\sigma\tilde\delta+2} &\displaystyle  \quad\text{if}\quad  \left(\tau\leq\tau^* \text{ and } \frac{\tau\tilde\gamma}{\tilde\delta}<\sigma\leq\frac{1}{L_K^2\tau}+\frac{\tilde\gamma}{L_K^2}\right).
\end{cases}
\end{equation}
In both cases, the minimum is reached when~$\sigma$ is maximal, equal to its upper bound given by~\eqref{eq:conditionsigma}. Hence,
\begin{equation}
\min_{\substack{\sigma\text{ subject to }\eqref{eq:conditionsigma}\\\theta\text{ subject to }\eqref{eq:conditionstep}}}\frac{\theta+2}{\sigma+1} = 
\begin{cases}
\displaystyle 1-\frac{\tilde\delta}{L_K^2\tau} & \!\!\!\text{if } \tau>\tau^*\\[3mm]
\displaystyle \min\left\{\frac{1}{\tau\tilde\gamma+1} ,\frac{1}{\tau\tilde\gamma+1}\,\frac{\tau\tilde\gamma+2}{\tilde\delta/(L_K^2\tau)+\tilde\delta\tilde\gamma/L_K^2+2}\right\}& \!\!\!\text{if } \tau\leq\tau^*.
\end{cases}
\end{equation}

Compare it to~$1/(\tau\tilde\gamma+1)$, and deduce the lower bound~$\omega^*(\tau)$:
\begin{equation}
\omega^*(\tau)
 = \max\left\{\frac{1}{\tau\tilde\gamma+1},\min_{\substack{\sigma\text{ subject to }\eqref{eq:conditionsigma}\\\theta\text{ subject to }\eqref{eq:conditionstep}}}\left\{\frac{\theta+2}{\sigma+1}\right\}\right\}.
\end{equation}
Thanks to~\eqref{eq:remark1}, it follows that
\begin{equation}
\omega^*(\tau)
 = 
\begin{cases}
\displaystyle 1-\frac{\tilde\delta}{L_K^2\tau} & \quad\text{if}\quad \tau>\tau^*\\[3mm]
\displaystyle \frac{1}{\tau\tilde\gamma+1} & \quad\text{if}\quad \tau\leq\tau^*.
\end{cases}
\end{equation}
In the second case,~$\sigma\tilde\delta$ is supposed to be greater than~$\tau\tilde\gamma$, so~$(\theta^*+1)/(\sigma\tilde\delta+2)$ is always smaller than~$1/(\tau\tilde\gamma+1)$. Therefore, the best rate is bounded from below by~$1/(\tau\tilde\gamma+1)$. Eventually, we get the following best rate:
\begin{equation}
\omega^*(\tau) = 
\begin{cases}
\displaystyle 1-\frac{\tilde\delta}{L_K^2\tau} & \quad\text{if}\quad \tau>\tau^*\\[3mm]
\displaystyle \frac{1}{\tau\tilde\gamma+1} & \quad\text{if}\quad \tau\leq\tau^*
\end{cases}
\end{equation}
which is minimal for~$\tau=\tau^*$.
This eventually leads to the best rate
\begin{equation}
\omega^* = 
1-\frac{\tilde\delta}{L_K^2\tau^*} =
\frac{1}{\tau^*\tilde\gamma+1} 
 = \frac{\sqrt{1+(4\,L_K^2)/(\tilde\gamma\tilde\delta)}-1
}{\sqrt{1+(4\,L_K^2)/(\tilde\gamma\tilde\delta)}+1}
 = \frac{\sqrt{1+4\kappa_F}-1
}{\sqrt{1+4\kappa_F}+1}
\end{equation}
obtained when~$\tau=\tau^*$ and~$\sigma=\tau^*\tilde\gamma/\tilde\delta$. $\square$

This choice leads to the following value for the solution error rate $\tilde\omega$:
\begin{equation}
\tilde\omega = 
\frac{1}{2\tau^*\tilde\gamma+1} 
 = \frac{\sqrt{1+(4\,L_K^2)/(\tilde\gamma\tilde\delta)}-1
}{\sqrt{1+(4\,L_K^2)/(\tilde\gamma\tilde\delta)}+3}
 = \frac{\sqrt{1+4\kappa_F}-1
}{\sqrt{1+4\kappa_F}+3}.
\end{equation}

Once again, the same computations prove that the best solution error rate $\tilde\omega$ is reached when 
\begin{equation}
\tilde\tau= \frac{\tilde\delta}{L_K^2}\left(1+\sqrt{1+\frac{L_K^2}{\tilde\gamma\tilde\delta}}\right)
\qquad\text{and}\qquad
\tilde\sigma= \frac{\tilde\gamma}{L_K^2}\left(1+\sqrt{1+\frac{L_K^2}{\tilde\gamma\tilde\delta}}\right)
\end{equation}
and leads to
\begin{equation}
\tilde\omega^* = \theta = 
 \frac{\sqrt{1+L_K^2/(\tilde\gamma\tilde\delta)}-1
}{\sqrt{1+L_K^2/(\tilde\gamma\tilde\delta)}+1}=
 \frac{\sqrt{1+\kappa_F}-1
}{\sqrt{1+\kappa_F}+1}
< \omega^*.
\end{equation}

\subsection{Overrelaxation on the dual variable}

Thanks to the symetry of Problem \eqref{eq:Pprimaldual}, similar results still hold if the relaxation is done on the dual variable~$y$ instead of the primal variable~$\xi$, namely if the updates are replaced by
\nopagebreak
\begin{equation}
\begin{cases}
\xi_{n+1} = \text{prox}_{\tau G}(\xi_n - \tau K^*\bar y_n)\\[2mm]
y_{n+1} = \text{prox}_{\sigma H^*}(y_n + \sigma K \xi_{n+1})\\[2mm]
\bar y_{n+1} = y_{n+1} + \theta\,(y_{n+1} - y_n).
\end{cases}
\end{equation}
As seen in~\eqref{eq:algoADMM-PD}, such an overrelaxation will be useful for the analysis of the ADMM.
It is equivalent to inverting the role of the dual and the primal variables. Indeed, Problem \eqref{eq:Pprimaldual} can be rewritten
\begin{equation}
\min_{y\in Y}\sup_{\xi\in X} \Big\{ H^*(y)-\langle K^*y,\xi\rangle-G(\xi)\Big\}
\end{equation}
which shares the same regularity assumptions as Problem \eqref{eq:Pprimaldual}. Hence, applying Theorem \ref{thm1} yields the following result:

\begin{theorem} \label{thm:relaxY}
Assume problem~\eqref{eq:Pprimaldual} has a solution, which is a saddle-point of~$\mathcal{L}$, denoted by~$(\xi^*,y^*)$.
Choose~$\tau>0$,~$\sigma>0$ and~$0<\theta\leq 1$ such that
\begin{equation}
\max\left\{\frac{1}{\tau\tilde\gamma+1},\frac{1}{\sigma\tilde\delta+1}\right\}\leq\theta\leq\frac{1}{L_K^2\tau\sigma}.
\label{eq:conditionstepY}
\end{equation}
Then, for any~$\omega$ such that
\begin{equation}
\max\left\{\frac{\theta+1}{\tau\tilde\gamma+2},\frac{1}{\sigma\tilde\delta+1}\right\}\leq\omega\leq\theta
\label{eq:conditionomegaY0}
\end{equation}
we have the following majoration for any~$N\in\mathbb{N}$ and any~$(\xi,y)\in \Xi\times Y$:
\begin{equation}
\begin{split}
(1-\omega L_K^2\tau\sigma)\frac{1}{2\tau}\,\lVert \xi_N - \xi \rVert^2
&+\frac{1}{2\sigma}\,\lVert y_N - y\rVert^2
\\&+\sum_{n=1}^N\frac{\omega^n}{\omega^{n-1}}\big(\mathcal{L}(\xi_n;y) - \mathcal{L}(\xi;y_n)\big)
\\[2mm]&\leq
\frac{\omega^N}{2\tau}\,\lVert \xi_0 - \xi \rVert^2 + \frac{\omega^N}{2\sigma}\,\lVert y_0 - y \rVert^2.
\end{split}
\label{eq:CVgapY}
\end{equation}
Now, define 
\begin{equation}
T_N := \sum_{n=1}^N \frac{1}{\omega^{n-1}}=\frac{1-\omega^N}{\omega^{N-1}(1-\omega)}
\end{equation}
and let
\begin{equation}
\Xi_N := \frac{1}{T_N}\sum_{n=1}^N\frac{1}{\omega^{n-1}} \, \xi_n
\qquad\text{and}\qquad
Y_N := \frac{1}{T_N}\sum_{n=1}^N\frac{1}{\omega^{n-1}} \, y_n.
\end{equation}
Then we have the following bound for any~$(\xi,y)\in \Xi\times Y$:
\begin{equation}
\begin{split}
\frac{1-\omega}{\omega(1-\omega^N)}\,(1-\omega L_K^2\tau\sigma)\frac{1}{2\tau}\,\lVert \xi-\xi_N \rVert^2&
+\frac{1-\omega}{\omega(1-\omega^N)}\,\frac{1}{2\sigma}\,\lVert y-y_N\rVert^2\\[2mm]&
+\mathcal{L}(\Xi_N;y) - \mathcal{L}(\xi;Y_N)\\[2mm]&
\leq
\frac{1}{T_N}\,\frac{1}{2\tau}\,\lVert \xi-\xi_0 \rVert^2
 + \frac{1}{T_N}\,\frac{1}{2\sigma}\,\lVert y-y_0 \rVert^2.
\end{split}
\label{eq:boundY}
\end{equation}
\end{theorem}
Note that the conditions on the parameters now slightly differ from the previous case. A variant can be found in~\cite[Appendix C2]{chambolle2016introduction}.

As in the previous case where the overrelaxation is done over the primal variable, we can prove the following result for the linear convergence of the solution errors:
\begin{corollary}\label{lemma2newY}
Assume problem~\eqref{eq:Pprimaldual} has a solution, which is a saddle-point of~$\mathcal{L}$, denoted by~$(\xi^*,y^*)$.
Let $(\xi_n,y_n)_n$ be generated by Algorithm~\eqref{eq:algoADMM-PD22}. Suppose there exist $\tau$, $\sigma$, $\theta$ and $\omega$ satisfying both conditions
\begin{equation}
\max\left\{\frac{1}{2\tau\tilde\gamma+1},\frac{1}{2\sigma\tilde\delta+1}\right\}\leq\theta\leq\frac{1}{L_K^2\tau\sigma}.
\label{eq:conditionstepY2}
\end{equation}
Then, for any~$\tilde\omega$ such that
\begin{equation}
\max\left\{\frac{\theta+1}{2\tau\tilde\gamma+1},\frac{1}{2\sigma\tilde\delta+2}\right\}\leq\tilde\omega\leq\theta.
\label{eq:conditionomegaY}
\end{equation}
Then, for any $N\in\mathbb{N}$, we have
\begin{equation}
\lVert y^*-y_N \rVert^2
\leq
\tilde\omega^N\left(\lVert y^*-y_0 \rVert^2 + \frac{\tau}{\sigma}\,\lVert \xi^*-\xi_0 \rVert^2\right).
\end{equation}
Moreover, if $\tilde\omega L_K^2\tau\sigma\neq 1$, then we also have
\begin{equation}
\lVert \xi^*-\xi_N\rVert^2
\leq
\frac{\tilde\omega^N}{1-\tilde\omega L_K^2\tau\sigma}\left(\frac{\sigma}{\tau}\lVert y^*-y_0 \rVert^2 + \lVert \xi^*-\xi_0 \rVert^2\right).
\end{equation}
\end{corollary}
Similar computations as in the previous section show that the best rate $\omega^*$ is achieved when choosing the following parameters:
\begin{equation}
\tau = \frac{\tilde\delta}{2\,L_K^2}\left(1+\sqrt{1+\frac{4\,L_K^2}{\tilde\gamma\tilde\delta}}\right)
\qquad\text{and}\qquad
\sigma = \frac{\tilde\gamma}{2\,L_K^2}\left(1+\sqrt{1+\frac{4\,L_K^2}{\tilde\gamma\tilde\delta}}\right)
\end{equation}
and, with $\kappa_F = L_K^2/(\tilde\gamma\tilde\delta)$, 
\begin{equation}
\theta = 
\frac{\displaystyle \sqrt{1+(4\,L_K^2)/(\tilde\gamma\tilde\delta)}-1}{\displaystyle \sqrt{1+(4\,L_K^2)/(\tilde\gamma\tilde\delta)}+1} 
= \frac{\displaystyle \sqrt{1+4\,\kappa_F}-1}{\displaystyle \sqrt{1+4\,\kappa_F}+1} <1
\end{equation}
which leads to $\omega^* = \theta$.

\section{Application : convergence rate for the ADMM in the smooth case}\label{sec:cvADMM}
As the ADMM is nothing but a particular instance of the oPDHG method with additional constraints on the parameter choice, its convergence rate is expected to be worse than that of the latter. In subsection~\ref{sub:uA}, it will indeed be derived from the computations of the previous section and shown to be greater than that of the oPDHG method.

However, as we will show it in the subsection~\ref{sub:aA}, it is possible to recover the same convergence rate as in the oPDHG method by introducing a slight modification in the ADMM iterations.

\subsection{Unaccelerated ADMM}\label{sub:uA}
As recalled in Section~\ref{sub:equivalence}, the ADMM iterations~\eqref{eq:algoADMM}, which aim at solving the primal problem
\begin{equation}
\min_{x\in X}  \Big\{f(x) := g(x)+h(Ax)\Big\}
\end{equation}
are equivalent to the oPDHG iterations~\eqref{eq:algoADMM-PD} applied to the primal-dual problem
\begin{equation}
\min_{\xi\in Y}\sup_{y\in Y} \Big\{ \mathcal{L}(\xi;y):=g_A(\xi) + \langle \xi,y \rangle - h^*(y) \Big\}.
\end{equation}
Hence, to study the convergence of the ADMM, one can either apply Theorem \ref{thm1} or Theorem \ref{thm:relaxY}, depending on the overrelaxation choice, with~$G = g_A$ and~$H = h$, and the identity operator~$K=\text{Id}$, of norm~$L_K=1$. We recall that~$Kx_n=\xi_n$ and that $g_A(Kx) = g(x)$ for any $x\in X$. The functions~$G$ and~$H$ are proved to be respectively~$\tilde\gamma=\gamma/L_A^2$-convex and~$\tilde\delta=\delta$-convex. 
In the case considered here, the relaxation is done on the dual variable, of parameter~$\theta = 1$.
The stepsize for the primal (resp. dual) proximal ascent is~$\tau>0$ (resp.~$\sigma=1/\tau$).

\subsubsection{Ergodic linear convergence}
Apply Theorem~\ref{thm:relaxY}. Parameters~$\theta$ and~$\sigma$ being constrained as stated above, Theorem~\ref{thm:relaxY} ensures that, provided one can find~$\tau>0$ such that
\begin{equation}
\max\left\{\frac{1}{\tau\gamma/L_A^2+1},\frac{1}{\delta/\tau+1}\right\}\leq1
\label{eq:conditionstepYuA}
\end{equation}
for any~$\omega$ such that
\begin{equation}
\max\left\{\frac{2}{\tau\gamma/L_A^2+2},\frac{1}{\delta/\tau+1}\right\}\leq\omega\leq1
\label{eq:conditionomegaYuA2}
\end{equation}
we have the following bound for any~$(x,y)\in X\times Y$:
\begin{equation}
\begin{split}
\frac{(1-\omega)^2}{\omega(1-\omega^N)}\frac{1}{2\tau}\,\lVert Ax-Ax_N \rVert^2&
+\frac{1-\omega}{\omega(1-\omega^N)}\,\frac{1}{2/\tau}\,\lVert y-y_N\rVert^2\\[2mm]&
+\mathcal{L}(AX_N;y) - \mathcal{L}(Ax;Y_N)\\[2mm]&
\leq
\frac{1}{T_N}\,\frac{1}{2\tau}\,\lVert Ax-Ax_0 \rVert^2
 + \frac{1}{T_N}\,\frac{1}{2/\tau}\,\lVert y-y_0 \rVert^2.
\end{split}
\label{eq:uA}
\end{equation}
We recall that $\mathcal{L}(Ax;y) = g(x)+\langle Ax,y\rangle-h^*(y)$ and $f(x) = \sup_{y\in Y}\mathcal{L}(Ax;y)$.
First note that, if we apply this inequality to $(x,y)=(x^*,y^*)$, then its left-hand side is nonnegative. Hence, the linear convergence of the dual iterates comes naturally. However, though the strong convexity ensures the convergence of the primal iterates~$x_N$, their convergence speed is not clear. We can solely estimate the convergence of~$Ax_N$, which is linear. Thanks to 
\begin{equation}
z_{N+1}-z^* = Ax_{N+1} - Ax^* + \tau\,(y^*-y_{N+1}) + \tau\,(y_N-y^*)
\end{equation}
we can nevertheless deduce the linear convergence of the primal iterates $z_N$. This also implies the linear convergence for the feasibility error $Ax_N-z_N$. If we now apply~\eqref{eq:uA} to $x=x^*$, using 
\begin{equation}
f(x^*) = \mathcal{L}(Ax^*;y^*)
 = \sup_{y\in Y}\mathcal{L}(Ax^*;y)
 \geq \mathcal{L}(Ax^*;Y_N)
\end{equation}
we get for any $y\in Y$
\begin{equation}
\begin{split}
\frac{(1-\omega)^2}{\omega(1-\omega^N)}\frac{1}{2\tau}\,\lVert Ax^*-Ax_N \rVert^2&
+\frac{1-\omega}{\omega(1-\omega^N)}\,\frac{1}{2/\tau}\,\lVert y-y_N\rVert^2\\[2mm]&
+\mathcal{L}(AX_N;y) - f(x^*)\\[2mm]&
\leq
\frac{1}{T_N}\,\frac{1}{2\tau}\,\lVert Ax^*-Ax_0 \rVert^2
 + \frac{1}{T_N}\,\frac{1}{2/\tau}\,\lVert y-y_0 \rVert^2.
\end{split}
\label{eq:uAstar}
\end{equation}
Let define $y^*_N\in Y$ as
\begin{equation}
y^*_N = \arg\max_{y\in Y}\mathcal{L}(AX_N;y)
\end{equation}
so that $\mathcal{L}(AX_N;y^*_N) = f(X_N)$. The left-hand side in \eqref{eq:uAstar} is then nonnegative for~$y=y^*_N$ and yields
\begin{equation}
0\leq f(X_N) - f(x^*)
\leq
\frac{1}{T_N}\,\frac{1}{2\tau}\,\lVert Ax^*-Ax_0 \rVert^2
 + \frac{1}{T_N}\,\frac{1}{2/\tau}\,\lVert y^*_N-y_0 \rVert^2.
\end{equation}
Hence, if the quantity $\lVert y^*_N-y_0\rVert$ is proved to be bounded, then the ergodic linear convergence of the ADMM in terms of objective error follows.
One can check that~$y^*_N = \nabla h(AX_N)$ and $y^*=\nabla h(Ax^*)$ which implies, thanks to the Lipschitz continuity of~$\nabla h$, that $\lVert y^*_N-y_0 \rVert\leq \lVert y^*-y^*_N \rVert+\lVert y^*-y_0 \rVert\leq \lVert Ax^*-AX_N \rVert/\delta+\lVert y^*-y_0 \rVert$. However, Equation \eqref{eq:uAstar} applied to $y=y^*$ and $N=n$ implies that
\begin{equation}
\lVert Ax^*-Ax_n \rVert^2
\leq
\frac{\omega(1-\omega^n)}{(1-\omega)^2}\,\frac{2\tau}{T_n}\left(\frac{1}{2\tau}\,\lVert Ax^*-Ax_0 \rVert^2
 + \frac{1}{2/\tau}\,\lVert y-y_0 \rVert^2\right).
\end{equation}
Thus, using the definition of $X_N$ and the convexity of the quadratic norm, we get
\begin{align}
\lVert Ax^*-AX_N \rVert^2
&= \left\lVert Ax^*-A\left(\frac{1}{T_N}\sum_{n=1}^N\frac{1}{\omega^{n-1}}\,x_n\right) \right\rVert^2\\
&\leq \frac{1}{T_N}\sum_{n=1}^N\frac{1}{\omega^{n-1}}\,\left\lVert Ax^*-Ax_n \right\rVert^2\\
\lVert Ax^*-AX_N \rVert^2
&\leq \frac{2\tau N}{T_N}\,\frac{\omega}{1-\omega}\left(\frac{1}{2\tau}\,\lVert Ax^*-Ax_0 \rVert^2
 + \frac{1}{2/\tau}\,\lVert y-y_0 \rVert^2\right)
\end{align}
thanks to the definition of $T_n$ \eqref{eq:defTN}. Since $(N/T_N)_N$ goes to zero, we get the desired result.

\subsubsection{Convergence rate}
Let us estimate the best convergence rate which can be achieved by the ADMM. 
Condition~\eqref{eq:conditionstepYuA} is always true.
Hence, for any~$\tau>0$, the convergence rate satisfies
\begin{equation}
\max\left\{\frac{1}{(\tau\gamma)/(2L_A^2)+1},\frac{1}{\delta/\tau+1}\right\}
\leq \omega\leq 1.
\end{equation}
The lower bound is equal to~$1/((\tau\gamma)/(2L_A^2)+1)$ when~$\tau\leq\sqrt{(2\delta L_A^2)/\gamma}$ and is equal to~$1/(\delta/\tau+1)$ otherwise. This leads to the best rate
\begin{equation}
\omega^*
 = \frac{1}{\sqrt{(\gamma\delta)/(2L_A^2)} +1}
 = \frac{1}{\sqrt{1/(2\kappa_f)} +1}
\quad\text{reached when}\quad
\tau = \sqrt{\frac{2\delta L_A^2}{\gamma}}.
\end{equation}
We call this parameter the \emph{optimal parameter} for the ADMM. Using this parameter also yields the following theoretical rate for the dual variable and $Ax_n$, given by Corollary~\ref{lemma2newY}:
\begin{equation}
\tilde\omega = 
\max\left\{\frac{1}{(\tau\gamma)/L_A^2+1},\frac{1}{2\delta/\tau+1}\right\}
 = \frac{1}{\sqrt{(2\gamma\delta)/L_A^2} +1}
 = \frac{1}{\sqrt{2/\kappa_f} +1}.
\end{equation}
This value can be easily proved to be the optimal one for $\tilde\omega$.

\subsection{Accelerated ADMM}\label{sub:aA}

We propose to relax the choice of step~$\tau$ in the updates of~$z$ and of~$y$ in the ADMM. Replacing~$\tau$ by~$\tau'\leq \tau$ in these two updates leads to the following algorithm:
\begin{equation}
\begin{cases}
\displaystyle x_{n+1} = \arg\min_{x\in X}\left\{g(x) + \langle Ax,y_n\rangle + \frac{1}{2\,\tau}\,\lVert Ax-z_n\rVert^2 \right\}\\[4mm]
\displaystyle z_{n+1} = \arg\min_{z\in Y}\left\{h(z) - \langle z,y_n\rangle + \frac{1}{2\,\tau'}\,\lVert Ax_{n+1}-z\rVert^2 \right\}\\[5mm]
\displaystyle y_{n+1} = y_n + \frac{1}{\tau'}\,(Kx_{n+1} - z_{n+1}).
\end{cases}
\label{eq:algoADMM2}
\end{equation}
\subsubsection{Equivalent oPDHG}
Following the same computations as in~\ref{sub:equivalence}, we show that iterations in Algorithm~\eqref{eq:algoADMM2} are equivalent to those of the following oPDHG algorithm
\begin{equation}
\begin{cases}
\xi^{n+1} = \text{prox}_{\tau g_A}\big(\xi^n - \tau\,\bar y^n\big)\\[2mm]
y^{n+1} = \text{prox}_{h^*/\tau'}(y^n+\xi^{n+1}/\tau')\\[2mm]
\bar y^{n+1} = y^{n+1}+\frac{\tau'}{\tau}\,(y^{n+1}-y^n)
\end{cases}
\end{equation}
where the relaxation parameter~$\theta = \tau'/\tau$ is linked to the ascent steps~$\tau$ and~$\sigma=1/\tau'$. 

Once again, Theorem~\ref{thm:relaxY} reads for any suitable~$\omega$,~$\tau$ and~$\tau'$:
\begin{equation}
\begin{split}
0\leq \frac{1-\omega}{\omega(1-\omega^N)}\,(1-\omega\tau/\tau')\frac{1}{2\tau}&\,\lVert Ax-Ax_N \rVert^2
+\frac{1-\omega}{\omega(1-\omega^N)}\,\frac{1}{2/\tau'}\,\lVert y-y_N\rVert^2\\[2mm]&
+\mathcal{L}(AX_N;y) - \mathcal{L}(Ax;Y_N)\\[2mm]&
\leq
\frac{1}{T_N}\,\frac{1}{2\tau}\,\lVert Ax-Ax_0 \rVert^2
 + \frac{1}{T_N}\,\frac{1}{2/\tau'}\,\lVert y-y_0 \rVert^2
\end{split}
\label{eq:aA}
\end{equation}
which yields a linear convergence in terms of objective error (in an ergodic sense). However, the best convergence rate achieved by the algorithm is expected to be better than that of the unaccelerated ADMM. Indeed, introducing the relaxed step $\tau'$ add a degree of freedom in the constraints over the value of $\omega$. Hence, it is minimized over a larger set and its minimal value is thus smaller.

Similarly to the unaccelerated case,~\eqref{eq:aA} ensures the linear convergence of the dual iterates. If $1-\omega\tau/\tau'$ does not cancel, it also implies the linear convergence of the primal iterates $(z_N)$. Otherwise, we lose the control on the convergence of $(Ax_N)$, thus on that of $(z_N)$.

\subsubsection{Convergence rate}
Let us derive the best convergence rate for Algorithm~\eqref{eq:algoADMM2}. We may use Theorem~\ref{thm:relaxY}, which ensures that steps~$\tau$ and~$\tau'$ are constrained by the relations
\begin{equation}
\max\left\{\frac{1}{\tau\gamma/L_A^2+1},\frac{1}{\delta/\tau'+1}\right\}\leq\frac{\tau'}{\tau}\leq1
\label{eq:conditions}
\end{equation}
and that the convergence rate is constrained by
\begin{equation}
\max\left\{\frac{\tau'/\tau+1}{\tau\gamma/L_A^2+2},\frac{1}{\delta/\tau'+1}\right\}\leq\omega\leq\frac{\tau'}{\tau}.\label{eq:conditionstautauprime}
\end{equation}
Hence, it is sufficient to find $(\tau,\tau')$ satisfying both~\eqref{eq:conditions} and~\eqref{eq:conditionstautauprime} which minimize the left-hand member in the latter.

One can also first use the remark made after Corollary~\ref{lemma2newY}. If no constraint on $\theta$ is made, then the best rate is achieved when 
\begin{equation}
\tau = \frac{\delta}{2}\left(1+\sqrt{1+\frac{4\,L_A^2}{\gamma\delta}}\right)
\qquad\text{and}\qquad
\sigma = \frac{1}{\tau'} = \frac{\gamma}{2\,L_A^2}\left(1+\sqrt{1+\frac{4\,L_A^2}{\gamma\delta}}\right)
\end{equation}
and
\begin{equation}
\theta = 
\frac{ \sqrt{1+4L_A^2/(\gamma\delta)}-1}{ \sqrt{1+4L_A^2/(\gamma\delta)}+1} .
\end{equation}
Let us check that such a choice satisfy $\theta = \tau'/\tau$. First, we have
\begin{equation}
\tau' = \frac{\displaystyle 2L_A^2/\gamma}{ 1+\sqrt{1+4L_A^2/(\gamma\delta)}}
 = \frac{ (\delta/2)(\sqrt{1+4L_A^2/(\gamma\delta)}+1)(\sqrt{1+4L_A^2/(\gamma\delta)}-1)}{ 1+\sqrt{1+4L_A^2/(\gamma\delta)}}
\end{equation}
which implies that
\begin{equation}
\frac{\tau'}{\tau}
 = \frac{ \sqrt{1+4L_A^2/(\gamma\delta)}-1}{ \sqrt{1+4L_A^2/(\gamma\delta)}+1}.
\end{equation}
Hence, these parameters can be chosen for the accelerated ADMM, and yields to the best rate. Thus, they are called \emph{optimal parameters} for the accelerated ADMM. 
With this parameter choice, we have~$\omega^* = \theta$. 
Note that the resulting rate is the same as the best one expected when applying the oPDHG on Problem~\eqref{eq:Pprimal0}. However, unlike in the oPDHG algorithm, this choice implies a loss of control on both $x$-iterates and~$z$-iterates. Moreover, this choice leads to the following rate $\tilde\omega$:
\begin{equation}
\tilde\omega = \frac{1}{2\delta/(\tau')^*+1}
=\frac{\sqrt{1+(4L_A^2)/(\gamma\delta)}-1}{\sqrt{1+(4L_A^2)/(\gamma\delta)}+3}
 = \frac{\sqrt{1+\kappa_f}-1}{\sqrt{1+\kappa_f}+3}
\end{equation}

To minimize the latter rate, we use the previous computations with $\gamma$ and $\delta$ doubled, which leads to the parameter choice
\begin{equation}
\tau' = \frac{\delta}{2}\left(\sqrt{1+\frac{L_A^2}{\gamma\delta}}-1\right)
\quad\text{and}\quad
\tau = \frac{\delta}{2}\left(\sqrt{1+\frac{L_A^2}{\gamma\delta}}+1\right)
\end{equation}
and the resulting rate:
\begin{equation}
\tilde\omega^*
 = \frac{\sqrt{1+L_A^2/(\gamma\delta)}-1}{\sqrt{1+L_A^2/(\gamma\delta)}+1}
 = \frac{\sqrt{1+\kappa_f}-1}{\sqrt{1+\kappa_f}+1}.
\end{equation}

\subsection{Theoretical rate comparison}
\begin{figure}\centering
\includegraphics[width = 0.7\linewidth]{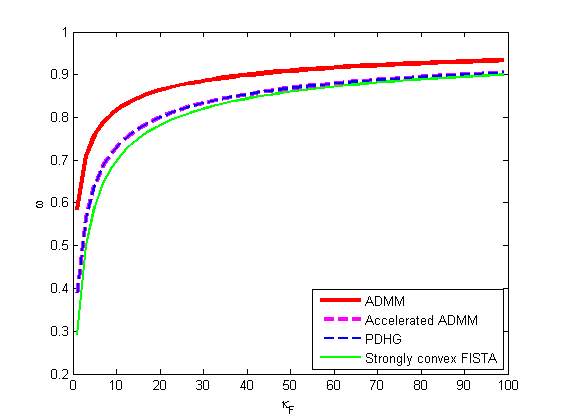}
\caption{theoretical rate comparison. In red/thick the unaccelerated ADMM, in magenta/thick dotted the accelerated ADMM, in blue/dotted the oPDHG and in green strongly convex FISTA with constant step.}
\label{fig:comptheo}
\end{figure}

Figure~\ref{fig:comptheo} compares the theoretical rates of the unaccelerated ADMM, the accelerated ADMM, the oPDHG method and strongly convex FISTA with constant step, by plotting for each algorithm the best rate with respect to the condition number $\kappa_f$. The rate achieved by strongly convex FISTA is the best one, but remains comparable with the accelerated ADMM and the oPDHG method. As expected, the unaccelerated ADMM yield larger rate values.

\section{Relations of other methods}\label{sec:relations}
In this section, we make a quick review on other linear convergence results for variant of the ADMM found in the literature. Generally, their differ from our result on the hypotheses made on the problem (both on the regularity of the objective function and on the operators).

\subsection{Overrelaxed ADMM}
In~\cite{nishihara2015general}, the authors propose to add an overrelaxation step in the spirit of \textsc{Nesterov}'s acceleration. They showed linear convergence rate when $h$ is assumed to be strongly convex and with \textsc{Lipschitz}-continuous gradient, while $B$ is invertible and $A$ is full column rank.

\subsection{Generalized ADMM} In~\cite{deng2016global}, the authors studied the ADMM in a wider framework, by allowing in each partial minimization to add an extra proximal term, which leads to a generalized ADMM.
Linear convergence rates are proved for four scenarios in which at least one of the functions $g$ or $h$ is strongly convex \emph{and} has a \textsc{Lipschitz} gradient, which is not assumed here. The case we treated is considered, but with extra assumptions (in particular, $h$ is supposed to be strongly convex). They provided an explicit convergence rate for only one scenario~\cite[Corollary 3.6]{deng2016global}.

\subsection{Relaxed ADMM}
It can be shown that the ADMM iterations are also equivalent to applying the \textsc{Douglas-Rachford} splitting (DRS) to the dual of~\eqref{eq:Pprimal0}. A relaxed version of the DRS, called \textsc{Peaceman-Rachford} splitting (PRS), can be obtained by introducing a relaxed parameter in the DRS iterations. Applying the PRS on the dual of~\eqref{eq:Pprimal0} hence leads to a so-called relaxed ADMM~\cite{davis2014faster}.
In~\cite[Theorem 6.3]{davis2014faster}, the authors proved the linear convergence rate of the relaxed ADMM in various cases (including the one we studied here), which depend on the assumptions made on the operators $A$ and $B$ (which is not supposed to be the negative identity) and / or on the regularity of the functions $g$ and $h$. However, the study is theoretical and does not provide explicit optimal rates.

\subsection{$K$-block ADMM}
In~\cite{hong2012linear}, the authors proved a linear convergence rate in the case where one can make assumptions on $g$ and $h$ which are supposed to be decomposable into a strictly convex term and a polyhedral one. This includes for instance the strongly convex case, but do not recover the smooth case studied in this paper. Furthermore, hypothesis on the rank of operators $A$ and $B$ (not necessary the negative identity) are made.
Moreover, their proof still holds when the objective function is a sum of $K$ separable convex functions (with an according number of variables).

\section{Applications}\label{sec:app}
\subsection{A toy example}
\subsubsection{Problem}
Let~$N$ be a integer. We consider the following constrained problem:
\begin{equation}
\min_{\substack{x=(x_i)_{i=0,\cdots,N-1}\in\mathbb{R}^N\\x_0=1}}\left\{
f(x):=\frac{M-m}{2}\,\lVert K_Nx\rVert_2^2 + \frac{m}{2}\,\lVert x\rVert^2_2
\right\}
\label{eq:5Pjouet}
\end{equation}
\noindent where the linear operator~$K_N:\mathbb{R}^N\to\mathbb{R}^{N-1}$ is defined by~$(K_Nx)_i=(x_{i+1}-x_{i})/2$ for any~$i=0,\cdots,N-2$, of norm~$\lVert K_N\rVert\leq 1$. The condition number of this problem is~$M/m$. Hence, if~$m$ is negligible compared to~$M$, then the problem is ill-conditioned. Let~$h(z):=(M-m)\,\lVert z\rVert^2_2/2$ for any~$z\in\mathbb{R}^{N-1}$ and~$g(x):=m\,\lVert x\rVert_2^2/2+\chi_{\{1\}}(x_0)$ for any~$x=(x_i)_{i=0,\cdots,N-1}\in\mathbb{R}^N$. The function $g$ is $m$-convex and~the convex conjugate $h^*:y\mapsto(M-m)^{-1}\lVert y\rVert_2^2/2$ is $(M-m)^{-1}$-convex.

\subsubsection{Solution}
\begin{figure}[t]
\centering
\includegraphics[width=0.7\linewidth]{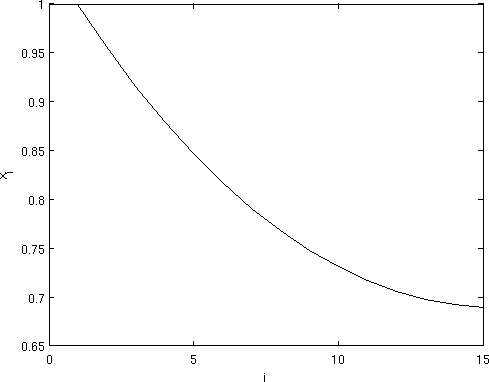}
\caption{Minimizer of~\eqref{eq:5Pjouet}.}
\label{fig:x}
\end{figure}

The minimizer of problem~\eqref{eq:5Pjouet} may be explicitly computed, by introducting the subvector~$\hat x$ given by:
\begin{equation}
\forall\,i=0,\cdots,N-2,\qquad
\hat x_i = x_{i+1}.
\end{equation}
\noindent such that~$x=(1,\hat x)$. The constrained problem~\eqref{eq:5Pjouet} can thus be rewritten in the unconstrained form
\begin{equation}
\min_{\hat x=(\hat x_i)_{i=0,\cdots,N-2}\in\mathbb{R}^{N-1}}\left\{
\frac{M-m}{2}\,\left(\lVert K_{N-1}\hat x\rVert_2^2+\frac{(\hat x_0-1)^2}{4}\right) + \frac{m}{2}\,(\lVert \hat x\rVert_2^2+1)
\right\}
\label{eq:5Pjouettilde}
\end{equation}
The minimizer~$\hat x^*$ is then given by the \textsc{Euler} equation, namely~$\hat x^* = A^{-1}b$ with
\begin{equation}
A = m\,\text{I}_{N-1} + (M-m)\, K_{N-1}^*K_{N-1} + \frac{M-m}{4}\,e_{0,0}
\end{equation}
\noindent where~$e_{0,0}$ denote the matrix of size~$N-1$ with null coefficients except the one at index~$(0,0)$ equal to 1. The vector $b$ is given by~$b:=(M-m)\,e_0/4$, with~$e_0$ the first vector of the canonical basis of~$\mathbb{R}^{N-1}$. Hence, the minimizer of the initial problem~\eqref{eq:5Pjouet} is~$x^* = (1,\hat x^*)$. For~$N=15$, $M=1000$, and~$m=1$, Figure~\ref{fig:x} plots~$x^*$.

\subsubsection{ADMM}
We apply the accelerated ADMM, which yields
\begin{equation*}
\begin{cases}
x_{n+1} = \displaystyle \arg\min_{\substack{x=(x_i)_i\in\mathbb{R}^N\\x_0=1}}\left\{\frac{m}{2}\,\lVert x\rVert_2^2 + \langle K_Nx,y_n\rangle + \frac{1}{2\,\tau}\,\lVert K_Nx-z_n\rVert_2^2 \right\}\\[4mm]
z_{n+1} = \displaystyle \arg\min_{z\in \mathbb{R}^{N-1}}\left\{\frac{M-m}{2}\,\lVert z\rVert_2^2 - \langle z,y_n\rangle + \frac{1}{2\,\tau'}\,\lVert K_Nx_{n+1}-z\rVert_2^2 \right\}\\[5mm]
y_{n+1} = \displaystyle y_n + \frac{1}{\tau'}\,(K_Nx_{n+1} - z_{n+1}).
\end{cases}
\end{equation*}
The~$z$-update is computed thanks to the \textsc{Euler} equation:
\begin{equation}
z_{n+1} = \frac{y_n+K_Nx_{n+1}/\tau'}{M-m+1/\tau'}.
\end{equation}
The~$x$-update is computed thanks to the subvectors we introduced above and is equivalent to minimizing
\begin{equation}
\frac{m}{2}\,\lVert \hat x\rVert_2^2 + \langle \hat x, K_{N-1}^*\hat y_n\rangle + \frac{1}{2\,\tau}\!\left(\lVert K_{N-1}\hat x-\hat z_n\rVert_2^2 +\left(\frac{\hat x_0-1}{2}-(z_n)_0\right)^{\!\!2}\right).
\end{equation}
The \textsc{Euler} equation ensures that~$\hat x_{n+1} = A_n^{-1}b_n$ with
\begin{equation}
A_n = m\,\text{I}_{N-1} + \frac{1}{\tau}\,K_{N-1}^*K_{N-1} + \frac{1}{4\,\tau}\,e_{0,0}
\end{equation}
\begin{equation}
b_n = - K_{N-1}^*\hat y_n + \frac{1}{\tau}\, K_{N-1}^*\hat z_n + \left(-\frac{(y_n)_0}{2}+\frac{1}{2\tau}\,(z_n)_0 + \frac{1}{4\tau}\right)e_0.
\end{equation}
We eventually have~$x_{n+1} = (1,\hat x_{n+1})$. 

\subsubsection{Parameters}
We tested two sets of parameters:
\begin{enumerate}
\item optimal parameter for the unaccelerated ADMM:
\begin{equation}\tau = \tau' = \sqrt{\frac{2}{m(M-m)}}
\end{equation}
(we assume that~$L=1$).
\item optimal paramaters for the accelerated ADMM: 
\begin{equation}
\tau = \frac{1}{2(M-m)}\left(\sqrt{1+\frac{4(M-m)}{m}}+1\right)
\end{equation}
\begin{equation}
\tau' = \tau-\frac{1}{M-m} = \frac{1}{2(M-m)}\left(\sqrt{1+\frac{4(M-m)}{m}}-1\right).
\end{equation}
\end{enumerate}
The convergence rates achieved in each case are respectively $1/(\sqrt{1/(M/m-1)/2}+1)$ and $(\sqrt{4M/m-3}-1)/(\sqrt{4M/m-3}+1)$.

\subsubsection{Comparison with oPDHG and strongly convex FISTA}
To solve problem~\eqref{eq:5Pjouettilde}, we can use the oPDHG iterations, by considering its primal-dual formulation
\begin{equation}
\min_{\substack{x=(x_i)_{i=0,\cdots,N-1}\in \mathbb{R}^N\\x_0=1}}\sup_{z'\in\mathbb{R}^{N-1}}\left\{
\frac{m}{2}\,\lVert x\rVert_2^2 + \langle Kx,z'\rangle - \frac{1}{2(M-m)}\,\lVert z'\rVert_2^2
\right\}.
\label{PD2}
\end{equation}
Hence, we are considering the following algorithm:
\begin{equation}
\begin{cases}
z'_{n+1} = \text{prox}_{\sigma h^*}(z'_n+\sigma\,K_N\bar x_n)\\[2mm]
x_{n+1} = \text{prox}_{\tau g}\big(x_n - \tau\,K_N^*z'_{n+1}\big)\\[2mm]
\bar x_{n+1} = x_{n+1}+\theta\,(x_{n+1}-x_n)
\end{cases}
\label{eq:algoADMM-PD2}
\end{equation}
for which the best theoretical convergence rate is achieved when choosing 
\begin{equation}
\tau = \frac{1}{2(M-m)}\left(1+\sqrt{1+\frac{4(M-m)}{m}}\right)
\end{equation}
\begin{equation}
\sigma = \frac{m}{2}\left(1+\sqrt{1+\frac{4(M-m)}{m}}\right)
\end{equation}

\begin{equation}
\theta = 
\frac{\displaystyle \sqrt{1+\frac{4(M-m)}{m}}-1}{\displaystyle \sqrt{1+\frac{4(M-m)}{m}}+1} <1.
\end{equation}
The $z'$-iterates are explicitly given by
\begin{equation}
z'_{n+1} = \frac{M-m}{M-m+\sigma}\,(z'_n+\sigma K_N\bar x_n)
\end{equation}
\begin{equation}
\hat x_{n+1} = \frac{\hat x_n/\tau -K_{N-1}^*\hat{z}'_{n+1}}{1/\tau+m}
\qquad\text{and}\quad
x_{n+1} = (1,\hat x_{n+1}).
\end{equation}
Note that, unlike in the ADMM iterations, there is no operator to invert.

We can also use the strongly convex FISTA algorithm, which solves problem~\eqref{eq:5Pjouet} by an accelerated FBS which can be written
\begin{equation}
\begin{cases}
x_{n+1} = \text{prox}_{\tau g}\big(\bar x_n-\tau\,\nabla h(\bar x_n)\big)\\[3mm]
\bar x_{n+1} = x_{n+1} + \theta_{n+1}\,(x_{n+1} - x_n).
\end{cases}
\end{equation}
where $\theta_n$ is given by~\eqref{eq:FISTApasconstant}, with $\tau=1/(M-m)$. The $x$-iterates are explicitly given by
\begin{equation}
\hat x_{n+1} = \frac{\widehat{\bar x_n}/\tau - (M-m)\widehat{(K^*_NK_N \bar x_n)}}{1/\tau+m}.
\end{equation}

\subsubsection{Results}
\begin{figure}
\centering
\includegraphics[width=0.95\linewidth]{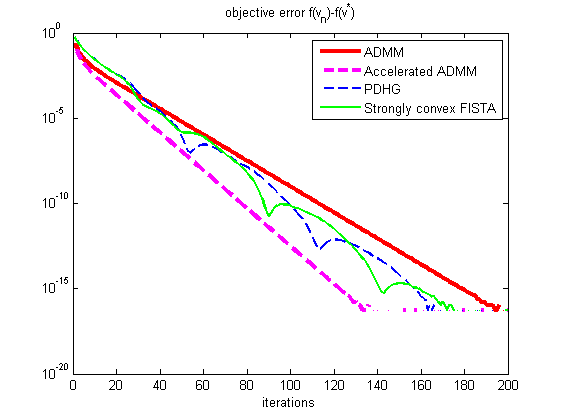}

(a) Objective error

\includegraphics[width=0.95\linewidth]{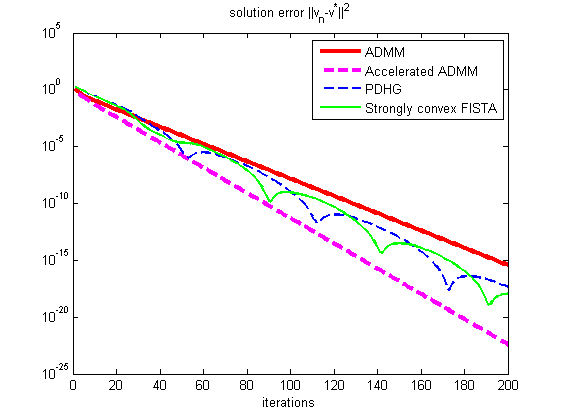}

(b) Solution error

\caption{Empirical convergence for the toy example.}\label{fig:Pjouet}
\end{figure}

To compare the convergence of each set of parameters, we used two tools:
\begin{enumerate}
\item the solution error~$\lVert x_n - x^*\rVert_2^2$;
\item the objective error~$f(x_n) - f(x^*)$.
\end{enumerate}
Figure~\ref{fig:Pjouet} displays the evolution of both measures, as well as the theoretical convergence decays expected in each case ($\tilde\omega$ and $\omega$). We chose $m=0.1$ and $M=10$, so that~$\kappa_f = 100$.

We first observe that, as expected, the accelerated ADMM has a better convergence than the unaccelerated ADMM. The empirical rates are better than the theoretical ones, which can be explained by the over-smoothness of the quadratic problem, compared to the assumptions required by the smooth case.

We also observe oscillations for both the oPDHG and strongly convex FISTA. Rippling for FISTA has been already observed for quadratic problems of this kind~\cite{o2015adaptive}. This phenomena occurs when the overrelaxation parameter $\theta$ is chosen too large compared to the eigenvalues of $m\,\text{I}_N+(M-m)K_N^*K_N$. Similar cause may explain the oscillations in the oPDHG, namely using overrelaxation steps can introduce oscillations when the according parameter are unproperly chosen. Hence, we do not expect to observe such oscillations for ADMM-like schemes.

\subsection{Denoising with TV-Huber}
\subsubsection{Problem}
\begin{figure}
\centering
\includegraphics[width = 0.3\linewidth]{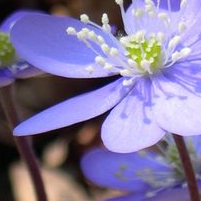}~
\includegraphics[width = 0.3\linewidth]{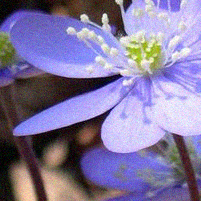}~
\includegraphics[width = 0.3\linewidth]{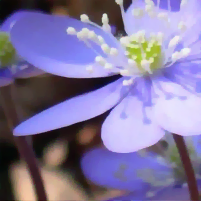}
\caption{We added a white Gaussian noise to an ideal image (left) to get noisy image (middle). The noise is of standard variation $10$ (the image values are between $0$ and $255$). The denoising is made by solving~\eqref{eq:5TVL2} (right).
Source : \emph{Hepatica nobilis flowers}, by Archenzo (detail).}
\label{fig:debruitageL2}
\end{figure}

We now apply the accelerated ADMM to a denoising problem, which is less smooth and more realistic than the toy example. Let~$g\in\mathbb{R}^{3N_xN_y}$ be a RGB-color (noisy) image. We want to solve the following problem:
\begin{equation}
\min_{v\in\mathbb{R}^{3N_xN_y}}\left\{
f(v):=\frac{\mu}{2}\,\lVert v-u\rVert_2^2 + h(\nabla v)
\right\}
\label{eq:5TVL2}
\end{equation}
where the gradient linear operator~$\nabla:\mathbb{R}^{3N_xN_y}\to\mathbb{R}^{3N_xN_y}\times\mathbb{R}^{3N_xN_y}$ is defined for any color image~$v$ by a pair of color images~$\nabla v = (\delta_xv,\delta_yv)^{\text T}$. The finite differences are given at any index~$(i,j)\in [0,N_x-1]\times[0,N_y-1]$ by
\begin{equation}
(\delta_xv)_{i,j}
= \begin{cases}
v_{i+1,j}-v_{i,j} & \text{if } i<N_x-1\\
0 & \text{otherwise}
\end{cases}
\end{equation}
and
\begin{equation}
(\delta_yv)_{i,j}
= \begin{cases}
v_{i,j+1}-v_{i,j} & \text{if } j<N_y-1\\
0 & \text{otherwise}.
\end{cases}
\end{equation}
The TV-\textsc{Huber} regularization term is defined by
\begin{equation}
h(\nabla v) = \sum_{i=0}^{N_x-1}\sum_{j=0}^{N_y-1} h_0\big(\lVert(\nabla v)_{i,j}\rVert\big)
\end{equation}
with
\begin{equation}
h_0(z) = \begin{cases}
\lvert z \rvert^2/2 & \text{if } \lvert z\rvert\leq 1\\
\lvert z \rvert-1/2 & \text{if } \lvert z\rvert> 1
\end{cases}\qquad\text{and}\qquad
h'_0(z) = \begin{cases}
z  & \text{if } \lvert z\rvert\leq 1\\
z/\lvert z \rvert & \text{if } \lvert z\rvert> 1.
\end{cases}
\end{equation} 
Hence, this term acts like a quadratic regularization when the image variations are small and like a TV regularization when they are larger (see Figure~\ref{fig:debruitageL2}). The quantity $\mu>0$ is a weight parameter.

The convex conjugate $h^*$ of the regularization function $h$ can be proved to be
\begin{equation}
h^*(y) = \sum_{i=0}^{N_x-1}\sum_{j=0}^{N_y-1} \left(\frac{1}{2}\,\lvert y_{i,j}\rvert^2 + \chi_{[0,1]}(\lvert y_{i,j}\rvert)\right)
\end{equation}
where $\chi_{[0,1]}(t)=0$ if $t\in[0,1]$ and $+\infty$ otherwise. This implies that $h^*$ is $1$-convex.

\subsubsection{ADMM}
Let~$g:=\mu\,\lVert\cdot-u\rVert_2^2/2$. We apply the accelerated ADMM to problem~\eqref{eq:5TVL2}, which leads to the following iterations:
\begin{equation*}
\begin{cases}
v_{n+1} = \displaystyle \arg\min_{v\in\mathbb{R}^{3N_xN_y}}\left\{
\frac{\mu}{2}\,\lVert v-u\rVert_2^2 + \langle \nabla v,\xi_n\rangle + \frac{1}{2\,\tau}\,\lVert \nabla v-\phi_n\rVert_2^2 
\right\}\\[4mm]
\phi_{n+1} = \displaystyle \arg\min_{\phi\in(\mathbb{R}^{3N_xN_y})^2}\left\{
h(\phi) - \langle \phi,\xi_n\rangle + \frac{1}{2\,\tau'}\,\lVert \nabla v_{n+1}-\phi\rVert_2^2 
\right\}\\[5mm]
\xi_{n+1} = \displaystyle \xi_n + \frac{1}{\tau'}\,(\nabla v_{n+1} - \phi_{n+1}).
\end{cases}
\end{equation*}
Each minimization is solved thanks to the \textsc{Euler} equation: the $v$-update reads
\begin{equation}
v_{n+1}
=\left(\mu\,\text{I} + \frac{1}{\tau}\,\nabla^*\nabla\right)^{-1}\left(\mu\,u+\frac{1}{\tau}\,\nabla^*\phi_n-\nabla^*\xi_n\right)
\end{equation}
whereas the $\phi$-update is given by
\begin{equation}
(\phi_{n+1})_{i,j}
= \frac{\tau'(\xi_n)_{i,j}+(\nabla v_{n+1})_{i,j}}{\lvert\tau'(\xi_n)_{i,j}+(\nabla v_{n+1})_{i,j}\rvert}\,\lvert(\phi_{n+1})_{i,j}\rvert
\end{equation}
with
\begin{equation}
\lvert(\phi_{n+1})_{i,j}\rvert
= \begin{cases}
\displaystyle \frac{\tau'\lvert(\xi_n)_{i,j}+(\nabla v_{n+1})_{i,j}\rvert}{\tau'+1}
 & \text{if }\lvert\tau'(\xi_n)_{i,j}+(\nabla v_{n+1})_{i,j}\rvert\leq \tau'+1 \\[3mm]
\displaystyle \lvert\tau'(\xi_n)_{i,j}+(\nabla v_{n+1})_{i,j}\rvert-\tau' 
 & \text{if }\lvert\tau'(\xi_n)_{i,j}+(\nabla v_{n+1})_{i,j}\rvert> \tau'+1.
\end{cases}
\end{equation}

\subsubsection{Parameters}
Before choosing the parameters, we recall the regularity of the problem. Functions~$h^*$ and~$g$ are respectively~$1$-convex and~$\mu$-convex. The gradient operator is bounded, of norm~$L\leq 2\sqrt{2}$ (this bound being tight when~$N_x$ or~$N_y$ go to~$+\infty$). Thus, we set $L=2\sqrt{2}$.
We tested two sets of parameters:
\begin{enumerate}
\item optimal parameter for the unaccelerated ADMM:~$\tau = \tau' = 4/\sqrt{\mu}$;
\item optimal paramters for the accelerated ADMM:
\begin{equation}
(\tau,\tau') = \left(\frac{1}{2}\left(\sqrt{1+\frac{32}{\mu}}+1\right),\frac{1}{2}\left(\sqrt{1+\frac{32}{\mu}}-1\right)\right).
\end{equation}
\end{enumerate}
These choices lead to the convergence rates $1/(\sqrt{\mu/8}+1)$ for the unaccelerated ADMM and $(\sqrt{1+32/\mu}-1)/(\sqrt{1+32/\mu}+1)$ for the accelerated one.
\subsubsection{oPDHG and strongly convex FISTA}
The primal-dual formulation of problem~\eqref{eq:5TVL2} is given by
\begin{equation}
\min_{v\in\mathbb{R}^{3N_xN_y}}\sup_{\phi\in(\mathbb{R}^{3N_xN_y})^2}\left\{
\frac{\mu}{2}\,\lVert v-u\rVert_2^2 + \langle\nabla v,\phi\rangle - h^*(\phi)
\right\}.
\label{TVL2PD}
\end{equation}
Hence using the oPDHG algorithm to solve it leads to the following iterations:
\begin{equation}
\begin{cases}
\phi'_{n+1} = \displaystyle \arg\min_{\phi'\in(\mathbb{R}^{3N_xN_y})^2}\left\{
h^*(\phi') + \frac{1}{2\sigma}\,\lVert \phi' - \phi'_n-\sigma\,\nabla\bar v_n\rVert^2
\right\}
\\[2mm]
v_{n+1} = \displaystyle \arg\min_{v\in\mathbb{R}^{3N_xN_y}}\left\{
\frac{\mu}{2}\,\lVert v-u\rVert^2 + \frac{1}{2\tau}\,\lVert v - v_n + \tau\,\nabla^*\phi'_{n+1}\rVert^2
\right\}\\[4mm]
\bar v_{n+1} =\displaystyle  v_{n+1}+\theta\,(v_{n+1}-v_n)
\end{cases}
\label{eq:algoADMM-PD23}
\end{equation}
which are computed thanks to the \textsc{Euler} equation:
\begin{equation}
(\phi'_{n+1})_{i,j} = \text{proj}_{[-1,1]}\!\left(\frac{(\phi'_n)_{i,j}+\sigma\,(\nabla\bar v_n)_{i,j}}{1+\sigma}\right)
\quad\text{and}\quad
v_{n+1} = \frac{v_n/\tau+\mu\,u - \nabla^*\phi'_{n+1}}{1/\tau+\mu}.
\label{eq:algoADMM-PD3}
\end{equation}
The best choice of parameters for this algorithm is (Theorem~\ref{thm2}):
\begin{equation}
\tau = \frac{1+\sqrt{1+32/\mu}}{16}
,~\quad
\sigma = \frac{1+\sqrt{1+32/\mu}}{16/\mu}
\quad\text{and}\quad
\theta = \frac{\sqrt{1+32/\mu}-1}{\sqrt{1+32/\mu}+1}.
\end{equation}

If we apply strongly convex FISTA to this problem, this leads to the following updates:
\begin{equation}
\begin{cases}
v_{n+1} = \displaystyle \arg\min_{v\in\mathbb{R}^{3N_xN_y}}\left\{
\frac{\mu}{2}\,\lVert v-u\rVert^2 + \frac{1}{2\tau}\,\lVert v - \bar v_n + \tau\,\nabla^*\nabla(\nabla h(\bar v_n))\rVert^2
\right\}\\[4mm]
\bar v_{n+1} = \displaystyle v_{n+1}+\theta_n\,(v_{n+1}-v_n)
\end{cases}
\end{equation}
which leads to the explicit update
\begin{equation}
v_{n+1} = \frac{\bar v_n/\tau + \mu\,u - \nabla^*\nabla(\nabla h(\bar v_n))}{1/\tau+\mu}.
\end{equation}
The variable relaxation parameter follows the update rule~\eqref{eq:FISTApasconstant} with $\tau=1/8$.

\subsubsection{Results}
To measure the convergence of the algorithm, we used the same two tools as in the previous case: the solution error and the objective error.

\begin{figure}
\centering
{\includegraphics[width=0.95\linewidth]{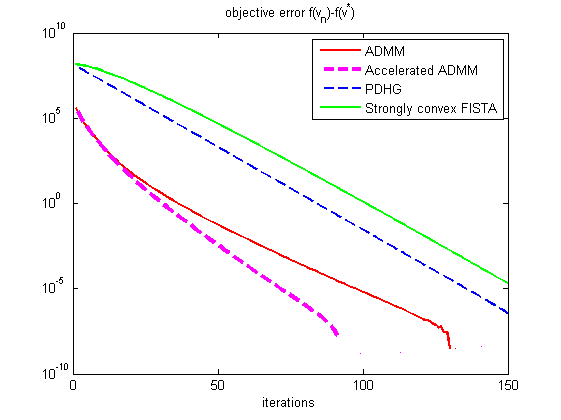}}

(a) Objective error

{\includegraphics[width=0.95\linewidth]{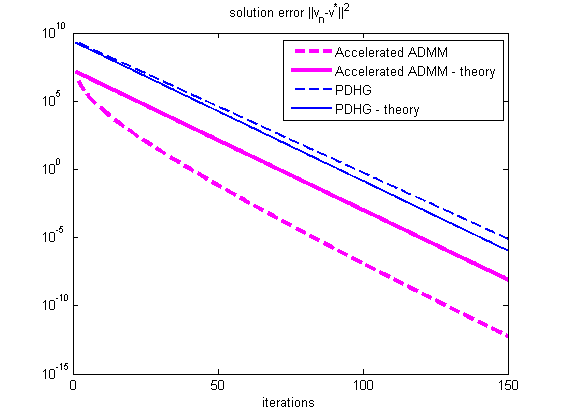}}

(b) Solution error - Comparison with the theoretical rate
\caption{Empirical convergence for TV-Huber denoising.}\label{fig:PTVHuber}
\end{figure}

Figure~\ref{fig:PTVHuber}(a) displays the evolution of the objective error, while Figure~\ref{fig:PTVHuber}(b) shows the decay of the solution error, for the accelerated ADMM and the oPDHG method. In the latter, the theoretical linear rate $\tilde\omega$ is also plotted for comparison. We chose~$\mu=10$. The solution error decreases as expected for all methods except strongly convex FISTA, for which we did not estimate a finer theoretical rate for the solution error. In practice, it seems that it converges with same rate as the oPDHG. Hence, in terms of solution error convergence, the accelerated ADMM provides the best empirical decay. For the objective error, the accelerated ADMM, the oPDHG method and strongly convex FISTA yield comparable decay rate. However, one should keep in mind that both the unaccelerated ADMM and the accelerated ADMM require an operator inversion, unlike the oPDHG method and strongly convex FISTA. Hence, even if comparable number of iterations are needed to achieve convergence, the ADMMs iterations are more time consuming than the other methods and should be used only when the inversion of the operator can be implemented efficiently. 

\section{Conclusion}
In this work, we studied the convergence of the oPDHG scheme in the case where the composite problem has a strongly convex part and a differentiable with a Lipschitz continuous gradient part. Using the equivalence between this algorithm and the ADMM, we provided a new convergence analysis of the latter. This analysis allowed us to introduce an accelerated variant of the ADMM by changing the augmented Lagrangian parameter, which is proved to have same convergence rate as the oPDHG method. Hence, we showed that in the smooth case, the choice of the ADMM parameter(s) can be crucial in terms of convergence rate. Experimental results confirmed this theoretical analysis. In particular, it has been observed that the accelerated ADMM does not introduce oscillations in some cases, unlike the oPDHG algorithm and strongly convex FISTA, which are known to be in practice more efficient than the ADMM-like scheme, since they require no operator inversion.

\bibliographystyle{plain}
\bibliography{biblio_optimisation}
\end{document}

%% file: preambuleGeneral.tex
\usepackage[english]{babel} 
\usepackage[utf8]{inputenc} 
\usepackage[T1]{fontenc} 
\usepackage{aeguill} 

\usepackage{tikz}
\usepackage{xcolor} 
\definecolor{couleur}{RGB}{0,62,92}
\definecolor{couleur2}{RGB}{202,117,96}

\usepackage{graphicx}
\usepackage{subfig}
\usepackage{epic}

%% file: preambuleMath.tex
\usepackage{mathrsfs}	
\usepackage{amsmath}
\usepackage{amsfonts}
\usepackage{amssymb}
\usepackage{dsfont} 
\usepackage{stmaryrd}
\everymath{\displaystyle}



%% file: preambuleMiseenpage.tex
\usepackage{geometry}
\usepackage{lscape} 

\author{Pauline Tan}
\newcommand{\auteur}{Pauline Tan}

\usepackage[colorlinks=true,breaklinks=true,linkcolor=couleur,citecolor=couleur,pdftitle={\titre},pdfauthor={\auteur}]{hyperref}

\usepackage{sectsty}
  \sectionfont{\color{couleur}{}}
  \subsectionfont{\color{couleur}{}}
  \subsubsectionfont{\color{couleur}{}}

\usepackage{enumitem}
  \setitemize[1]{label=$\bullet$}
  \setenumerate[1]{font=\color{couleur}}

\usepackage{framed}
\usepackage[framed]{ntheorem}
  \makeatletter
  \newtheoremstyle{theocouleur}{\item[\hskip\labelsep {\color{couleur}\theorem@headerfont##1\ ##2}\theorem@separator]\mbox{}}{\item[\hskip\labelsep \theorem@headerfont##1\ ##2\ (##3)\theorem@separator]\mbox{}}
  \makeatother

  \theoremstyle{theocouleur}
  \newshadedtheorem{theorem}{Theorem}
  \newshadedtheorem{lemma}{Lemma}
  \newshadedtheorem{corollary}{Corollary}
  \newshadedtheorem{proposition}{Proposition}
  \newshadedtheorem{definition}{Définition}

\newenvironment{montitre}{
	\centering
	\bigbreak
	\scshape\LARGE\color{couleur}}{\bigbreak\normalcolor\bigbreak}

\newcommand{\monref}[1]{(\hyperref[eq:#1]{#1})}

\newcommand{\monrefplusone}[2]{(\hyperref[eq:#1]{#2})}

\makeatletter
\newcommand{\specialnumber}[1]{%
  \def\tagform@##1{\maketag@@@{(\ignorespaces##1\unskip\@@italiccorr#1)}}%
}
\newcommand{\specialeqref}[2]{\begingroup
  \def\tagform@##1{\maketag@@@{(\ignorespaces##1\unskip\@@italiccorr#2)}}%
  \eqref{#1}\endgroup}
\makeatother